\documentclass[12pt]{article}
\usepackage{amsmath,amssymb,graphicx}
\usepackage[latin1]{inputenc}
\begin{document}

\title{\vspace{2cm}{\bf Computation of the multi-chord distribution of convex and concave polygons}}
\author{
Ricardo Garc\'\i a-Pelayo\thanks{%
E-mail: r.garcia-pelayo@upm.es} \\
\\
ETS de Ingenier\'ia Aeron\'autica \\
Plaza del Cardenal Cisneros, 3 \\
Universidad Polit\'ecnica de Madrid \\
Madrid 28040, Spain\\}
\date{}
\maketitle

%

\begin{abstract}
Analytical expressions for the distribution of the length of chords corresponding to the affine invariant measure on the set of chords are given for convex polygons. These analytical expressions are a computational improvement over other expressions published in 2009 and 2011. The correlation function of convex polygons can be computed from the results obtained in this work, because it is determined by the distribution of chords.

An analytical expression for the multi-chord distribution of the length of chords corresponding to the affine invariant measure on the set of chords is found for non convex polygons. In addition we give an algorithm to find this multi-chord distribution which, for many concave polygons, is computationally more efficient than the said analytical expression. The results also apply to non simply connected polygons.
\end{abstract}
\bigskip

\medskip

\bigskip

\section{Introduction}

The distribution of the length of chords of a plane figure is a topic that became popular in connection with Bertrand's paradox in 1907 \cite{Kendall1963, Clark2002}. The study of this paradox led mathematicians to realize that different problems lead to different measures on the set of chords. The relations between these measures on the set of chords of convex sets of the Euclidean space have been studied \cite{Kingman1965, Kingman1969, Piefke1978}. Among these measures there is only one (up to a factor), usually denoted by $\mu$, which is invariant under translations and rotations (see the first three chapters of \cite{Deltheil1919} or chapter 1 of \cite{Kendall1963} or \cite{Santalo1976}). The distribution of the length of chords corresponding to $\mu$ (and other measures) had been found analytically for the equilateral triangle \cite{Sulanke1961}, the general triangle \cite{Duma2009,Ciccariello2010}, the circle, the rectangle \cite{Coleman1969}, the ellipse \cite{Piefke1979}, the regular pentagon \cite{Aharonyan2005}, the regular hexagon \cite{Harutyunyan2007}, the trapeze rectangle \cite{Duma2011} and the isosceles trapezium \cite{Sorrenti2012b}, to name only plane figures. In the last years the distribution for regular polygons has been  found \cite{Harutyunyan2009,Baesel2014}. In 2009 Ciccariello gave an analytical expression for the correlation function of convex polygons \cite{Ciccariello2009} and showed that its second derivative is the chord length distribution (Appendix C of \cite{Ciccariello2009}. In 2011 the correlation function of convex polygons was found by H. S. Harutyunyan and V. K. Ohanyan \cite{Harutyunyan2011} differently. First they found the distribution for each pair of sides of a convex polygon following an approach based on Pleijel's identity (\cite{Ambartzumian1990}, page 156). Then they wrote the distribution for convex polygons as the sum of the distributions for each pair of sides. The work of these three authors, however, doesn't furnish an efficient tool to compute distributions of the length of chords. For example in section III of \cite{Ciccariello2009} and in section 5 of \cite{Harutyunyan2011} one can see that it takes them a lot of work to compute the distribution of the chords' length for an isosceles right triangle and a rhombus, respectively.

In section 2 of this article we find the distribution for triangles by direct geometrical arguments. The equilateral case was found in 1961 (\cite{Sulanke1961}, page 57) and the general case was found in 2010 by Ciccariello \cite{Ciccariello2010}. However, his final expressions are much more involved than the ones presented here. Ciccariello has checked on an acute and an obtuse triangle that both methods match \cite{PrivateCommunication}. In section 3 we find the distribution of the length of chords of two concurrent segments, as a step towards finding the distribution for a pair of sides in section 4. That the latter distribution can be derived from the former has been acknowledged earlier \cite{Boettcher2012}. In order to find the distribution of the length of chords for two concurrent segments, we have used the distribution for triangles as a basis, instead of Ciccariello's approach \cite{Ciccariello2009} or Pleijel's identity (\cite{Ambartzumian1990}, page 156). This and the use of the indicator function of an interval yields mathematical expressions which are computationally more friendly then the corresponding ones in \cite{Harutyunyan2011}. In section 5 we also write the distribution for convex polygons as the sum of the distributions for each pair of sides and, in addition, we prepare the analytical expressions for their input to be Cartesian coordinates. The graph of the distribution of the length of chords for a rhombus can be produced in less than a second on a personal computer by simply giving the coordinates of its vertices. In section 6 we compute the distribution of the chords' length of the Mallows and Clark dodecagon \cite{Mallows1970}, which has theoretical and practical interest (see the Introduction of \cite{Gille2009}).

While much is known about the distribution of random chords in convex bodies \cite{Ren2014}, the distribution of the length of chords in non convex figures has features which are absent in the convex case. In fact, one can tell whether a body is convex or not from its distribution of chords (\cite{Sulanke1961}, p. 61, proposition 4). This is also apparent when one tries to find the distribution by Monte Carlo methods \cite{Vlasov2011}. Despite these difficulties the distribution of the length of chords of several non convex figures has been computed or found analytically \cite{Gille2001a, Gille2001b, Mazzolo2003, Baesel2008, Gruy2008, Ciccariello2009, Barrilla2011, Sorrenti2012, Gruy2014}. A distribution for the length of chords can be found by taking the second derivative of the correlation function \cite{Ciccariello1981, Vlasov2007, Ciccariello2009}. This distribution, however, is a signed distribution for non convex bodies, that is, the probability may take negative values. This happens for chords for which one end lies inside the body and the other end lies outside the body. Considerations of such chords follows naturally from the interpretation of the second derivative of the correlation function. When we restrict ourselves to chords which lie fully within the body one may, for concave bodies, consider the sum of the lengths of the interceptions of a straight line by the body as the random variable, and this is the One-Chord Distribution (OCD) \cite{Gille2000}.

Here we take what from a human point of view (perhaps not from a mathematical point of view) is the most natural choice: we consider only chords which lie fully within the body and consider each segment independently from others which may lie on the same straight line. This is called the Multi-Chord Distribution (MCD) \cite{Gille2000}. Here we find an analytical expression for the MCD corresponding to the invariant measure $\mu$ for non convex polygons. We also give an algorithm to find the MCD which is numerically more efficient than the said analytical expression and use it to compute the MCD for the five-pointed star. The results also apply to non simply connected polygons.

In sections 2-6 we treat the case of convex polygons and in sections 7-14 the case of concave polygons. We provide a link where two Mathematica files with implementations of the procedures used in the convex and concave cases, respectively, can be downloaded. \cite{RGPN}.

The distribution of chords' length has applications in the analysis of images provided by a microscope \cite{Lang2001} and others (see, for example, the first paragraph of \cite{Gille2000}). The importance of the correlation function of polygons in small angle scattering has been the subject of recent research \cite{Gille2000,Ciccariello2009}. The correlation function of convex polygons can be computed from the results obtained in this work, because the distribution of chords determines the correlation function (\cite{Piefke1978}, for a clear exposition of this in English see Appendix A of \cite{Vlasov2007}).

\section{Distribution of the length of chords in a triangle}

In this work we shall use repeatedly these two facts: the affine invariant measure $\mu$ of a set of parallel straight lines is proportional to its cross section and the affine invariant measure $\mu$ of a sheaf of straight lines is proportional to the angle that they span. See, for example, references \cite{Kendall1963,Kingman1969,Santalo1976} for details. We prefer to find the distribution of measure rather than the distribution of probability to avoid dealing with normalization factors. It is better to normalize, if desired, only after all calculations are done.

Since we are going to find the densities of measure of many different geometrical figures or bodies, we shall always use the symbol $\rho$ for it, which will refer to different figures or bodies in different sections. Sometimes different densities are used in the same section. In that case either subscripts will be used or else they may be distinguish by their arguments. For example in section 4 $\rho(A, A', B, B'; \ell)$ is the value of the measure density of the length of the chords that join the segments $AA'$ and $BB'$ at $\ell$, whereas $\rho(OA', OB', \gamma; \ell)$ is the value of the measure density of the length of the chords that join the segments $OA'$ and $OB'$ (where $\angle{A'OB'} = \gamma$) at $\ell$. Still to avoid heavy notation, given two points $A$ and $B$, by $AB$ we may mean either a segment of extremes $A$ and $B$ or the length of the said segment.

{\bf{Lemma 2.1.}} {\it{The longest chord parallel to any given direction hits one of the vertices of the triangle.}}

{\bf{Proof.}} Since the three vertices of a triangle are not aligned, as we move along a direction perpendicular to the set of parallel straight lines, we meet the three vertices at different times. When the first vertex is met the length of the chord is zero. Then it increases until the second vertex is met, after which it decreases again to be of length zero at the third vertex. $\Box$
\medskip

{\bf{Lemma 2.2.}} {\it{The length of the chords of a triangle parallel to any given direction is uniformly distributed between zero and its maximum length.}}

{\bf{Proof.}} Consider the set of chords of a triangle parallel to any given direction. Their length changes linearly with the displacement as we move perpendicularly to that direction. $\Box$

{\bf{Definition 2.1}} We say that a chord of a given direction is ``dominated" by the vertex $A$ when chords having that direction reach their maximum length when they meet the vertex A. We denote this set by $D(A)$.
\medskip

{\bf{Lemma 2.3.}} {\it{The chords of a triangle belong to $D(A) \cup D(B) \cup D(C)$. For any two vertices $A$ and $B$, $\mu(D(A) \cap D(B)) = 0$.}}

{\bf{Proof.}} The first statement is obvious. As for the second one, there is only one direction in $D(A) \cap D(B)$. $\Box$
\medskip

We now proceed to find the distribution of measure for the length of the chords dominated by a vertex, and later we may add the contributions of each vertex.
\medskip

\begin{figure}[ht]
  \begin{center}
  \includegraphics[width=0.75 \textwidth]{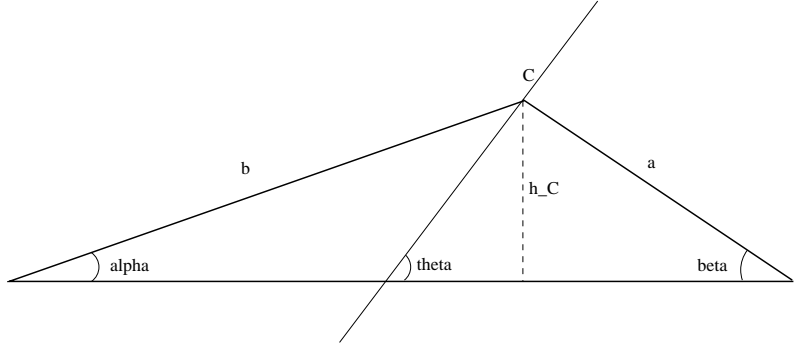}
  \caption{\label{Chords1}Vertices $A$ and $B$ are both acute}
  \end{center}
\end{figure}

When computing the distribution of measure of $D(C)$, vertices $A$ and $B$ may or may not be both acute. We consider these two cases separately.

In the first case, the set $D(C)$ is the set of chords whose direction $\theta \in [ \alpha, \pi - \beta ]$, as seen in the figure. We denote by $D(C, \theta)$ the set of chords of direction $\theta$ which are in $D(C)$, and by $\mu_C(\theta)$ its measure (restricted to the direction $\theta$). This measure is simply the cross section:

\begin{equation}\label{1}
\mu_C(\theta) = c \sin \theta.
\end{equation}
We note that
\begin{equation}
\mu(D(C)) = \int_\alpha^{\pi - \beta} d\theta\ c \sin \theta = c (\cos \alpha + \cos \beta),
\end{equation}
which, by Cauchy's theorem \cite{Cauchy1850, Czuber1884, Kendall1963, Santalo1976}, yields the identity
\begin{equation}
c (\cos \alpha + \cos \beta) + b (\cos \gamma + \cos \alpha) + a (\cos \beta + \cos \gamma) = a+b+c,
\end{equation}
which can be obtained by direct and simple means.

We denote by $h_C(\theta)$ the length of the longest chord in $D(C, \theta)$.
\begin{equation}\label{4}
h_C(\theta) = {h_C \over \sin \theta}.
\end{equation}

\begin{figure}[ht]
  \begin{center}
  \includegraphics[width=0.4 \textwidth]{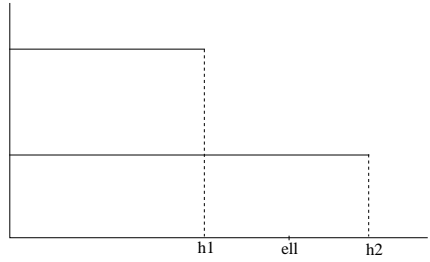}\caption{\label{Chords2}Density of measure for the length of the chords in $D(C, \theta_1)$ and $D(C, \theta_2)$}
  \end{center}
\end{figure}

The measure density of the length of the chords in $D(C, \theta)$ is constant when the length $\ell$ is smaller than $h_C(\theta)$, and zero otherwise. Since the measure of $D(C, \theta)$ is $\mu_C(\theta)$, the measure density (which we denote by $\rho_C$) is

\begin{equation}
\rho_C(\theta, \ell) = {\mu_C(\theta) \over h_C(\theta)} H(h_C(\theta) - \ell),
\end{equation}
where $H$ is the Heaviside function, which is 1 when its argument is positive and 0 otherwise. The measure density $\rho_C(\theta, \ell)$ contributes only to those $\ell$ which are smaller than $h_C(\theta)$. This is seen in Figure \ref{Chords2}, which depicts the measure densities for two different angles in $[ \alpha, \pi - \beta ]$ which have different $h_C(\theta)$. Only the measure density $\rho_C(\theta_2, \ell)$ contributes to the measure density at the $\ell$ shown in Figure \ref{Chords2}, which is between $h_C(\theta_1)$ and $h_C(\theta_2)$.

Formally, the distribution of measure for the length of the chords in $D(C)$ is

\begin{equation}\label{6}
\rho_C(\ell) = \int_{\alpha}^{\pi - \beta} d\theta\ {\mu_C(\theta) \over h_C(\theta)} H(h_C(\theta) - \ell).
\end{equation}
We are going to substitute the Heaviside function in the integrand by appropriate limits for the integral for each $\ell$. When $\ell < h_C$, all the $\rho_C(\theta, \ell)$ contribute to $\rho_C(\ell)$. When $h_C < \ell <\min(a,b)$, consider the two chords of length $\ell$ which meet at the vertex C. They form angles $\arcsin{h_C \over \ell}$ and $\pi - \arcsin{h_C \over \ell}$ with the side $c$, where by $\arcsin{h_C \over \ell}$ we mean the value between 0 and $\pi/2$. Then it is clear that the chords of length in $[h_C, \min(a,b)]$ belong only to the sets $D(C, \theta)$ s. t. $\theta \in [\alpha, \arcsin{h_C \over \ell}] \cup [\pi - \arcsin{h_C \over \ell},\pi - \beta]$. Likewise, if, say, $a < b$, and $a < \ell < b$, then the chords are in the sets $D(C, \theta)$ s. t. $\theta \in [\alpha, \arcsin{h_C \over \ell}]$. Then

$$
\rho_C(\ell) = H(h_C - \ell) \int_{\alpha}^{\pi - \beta} d\theta\ {\mu_C(\theta) \over h_C(\theta)} +
$$

\begin{equation}
H(\ell - h_C) H(b - \ell) \int_{\alpha}^{\arcsin{h_C \over \ell}} d\theta\ {\mu_C(\theta) \over h_C(\theta)} +
H(\ell - h_C) H(a - \ell) \int_{\pi - \arcsin{h_C \over \ell}}^{\pi - \beta} d\theta\ {\mu_C(\theta) \over h_C(\theta)}.
\end{equation}

To do the integral we use

\begin{equation}
h_C = {2 A \over c},
\end{equation}
where $A$ is the area of the triangle, formulae (\ref{1}) and (\ref{4}) which yield

\begin{equation}
{\mu_C(\theta) \over h_C(\theta)} = {c^2 \over 2 A} \sin^2 \theta
\end{equation}
and the definite integral

\begin{equation}
\int_{aa}^{bb} \sin^2 \theta = {1 \over 2} \Big( bb-aa+{\sin 2 aa - \sin 2 bb \over 2} \Big).
\end{equation}
The result is

$$
\rho_C(\ell) = {c^2 \over 4 A} \left[ H\left({2 A \over c} - \ell \right) \left( \gamma + {\sin 2 \alpha \over 2} + {\sin 2 \beta \over 2} \right) \right. +
$$

$$
H(\ell - {2 A \over c}) H(b - \ell) \left( \arcsin{2 A \over c \ell} - \alpha + {\sin 2 \alpha \over 2} - {2 A \over c \ell} \sqrt{1 - \left({2 A \over c \ell}\right)^2} \right)  +
$$

\begin{equation}\label{triagudo-2}
\left.
H(\ell - {2 A \over c}) H(a - \ell) \left( \arcsin{2 A \over c \ell} - \beta + {\sin 2 \beta \over 2} - {2 A \over c \ell} \sqrt{1 - \left({2 A \over c \ell}\right)^2} \right)
\right].
\end{equation}

Using the definition

\begin{equation}
\rho_{\rm{ac}}(\ell, \alpha, \beta, \gamma, a, b, c) \equiv \rho_C(\ell),
\end{equation}
we may write the measure density of length of chords of an acute triangle as

\begin{equation}\label{triagudo}
\rho_{\rm{ac.tr}}(\ell) = \rho_{\rm{ac}}(\ell, \beta, \gamma, \alpha, b, c, a) +
\rho_{\rm{ac}}(\ell,  \gamma, \alpha, \beta, c, a, b,) + \rho_{\rm{ac}}(\ell, \alpha, \beta, \gamma, a, b, c).
\end{equation}
We could, of course, have omitted the arguments $\alpha, \beta$ and $\gamma$, because they are are determined by $a, b$ and $c$, but we believe things are clearer this way.
\bigskip

\begin{figure}[ht]
  \begin{center}
  \includegraphics[width=0.55 \textwidth]{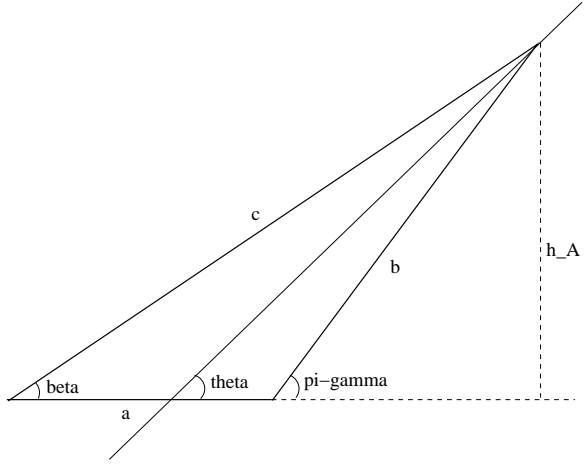}\caption{\label{Chords3}Vertices $B$ and $C$ are not both acute}
  \end{center}
\end{figure}

We now consider the case in which one of the vertices opposite to the vertex considered is obtuse. In Figure \ref{Chords3} the same obtuse triangle as in Figure \ref{Chords1} is considered, but now it is lying on side $a$. All reasonings leading to formula (\ref{6}) still hold. This formula for vertex $A$ is

\begin{equation}
\rho_A(\ell) = \int_{\beta}^{\pi - \gamma} d\theta\ {\mu_A(\theta) \over h_A(\theta)} H(h_A(\theta) - \ell).
\end{equation}
The reasonings after formula (\ref{6}) yield now only two terms, because the entire triangle lays now to only one side of the height $h_A$.

$$
\rho_A(\ell) = H(b - \ell) \int_{\beta}^{\pi - \gamma} d\theta\ {\mu_A(\theta) \over h_A(\theta)} +
H(\ell - b) H(c - \ell) \int_{\beta}^{\arcsin{h_A \over \ell}} d\theta\ {\mu_A(\theta) \over h_A(\theta)} =
$$

$$
{a^2 \over 4 A} \left[ H\left(b - \ell \right) \left( \alpha + {\sin 2 \beta \over 2} + {\sin 2 \gamma \over 2} \right) \right. +
$$

\begin{equation}\label{triobtuso-2}
\left.
H(\ell - b) H(c - \ell)
\left( \arcsin{2 A \over a \ell} - \beta + {\sin 2 \beta \over 2} - {2 A \over a \ell} \sqrt{1 - \left({2 A \over a \ell}\right)^2} \right)
\right].
\end{equation}

Using the definition

\begin{equation}
\rho_{\rm{obt}}(\ell, \alpha, \beta, \gamma, a, b, c) \equiv \rho_A(\ell),
\end{equation}
we may write the measure density of length of chords of an obtuse triangle as

\begin{equation}\label{triobtuso}
\rho_{\rm{obt.tr}}(\ell) = \rho_{\rm{obt}}(\ell, \alpha, \beta, \gamma, a, b, c) +
\rho_{\rm{obt}}(\ell, \beta, \alpha, \gamma, b, a, c) + \rho_{\rm{ac}}(\ell, \alpha, \beta, \gamma, a, b, c),
\end{equation}
where $c$ is the longest side.
\bigskip

We have proven

{\bf{Proposition 2.1.}} {\it{The density of measure of the length of chords of a triangle, when the set chords is endowed with the measure $\mu$, is given by expression (\ref{triobtuso}) when the triangle has an obtuse angle, and by expression (\ref{triagudo}) otherwise.}}

Recently Gasparyan and Ohanyan \cite{Gasparyan2013} have found the density of measure of the length of chords of a triangle. They have not, however, organized their results in a form as compact as here, where all that is needed to find the density of measure of the length of chords of a triangle is formulae (\ref{triagudo-2})-(\ref{triagudo}) and (\ref{triobtuso-2})-(\ref{triobtuso}).

\section{Distribution of the length of chords of two concurrent segments}

\begin{figure}[ht]
  \begin{center}
  \includegraphics[width=0.45 \textwidth]{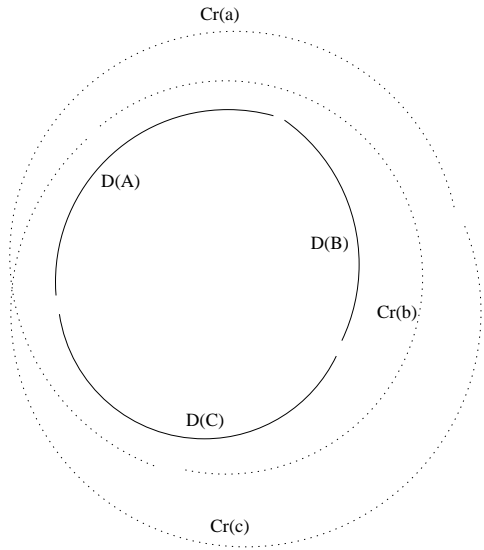}\caption{\label{Chords44}Relations among the sets $D$(vertex) and $Cr$(side)}
  \end{center}
\end{figure}

Let us denote by $Cr(a)$ the set of chords which cross the side $a$ of a triangle. From this definition these three statements follow: 1) $D(A) \subset Cr(a), D(B) \subset Cr(b)$ and $D(C) \subset Cr(c)$; 2) up to sets of measure zero, each chord is contained in exactly two of the sets $Cr(a)$, $Cr(b)$ and $Cr(c)$; 3) $D(C) = (D(C) \cap Cr(a)) \sqcup (D(C) \cap Cr(b)), D(B) = (D(B) \cap Cr(a)) \sqcup (D(B) \cap Cr(c))$ and $D(A) = (D(A) \cap Cr(b)) \sqcup (D(A) \cap Cr(c))$, where $\sqcup$ denotes disjoint (up to sets of measure zero) union. The relations among the sets $D$(vertex) and $Cr$(side), which were summarized at the beginning of these section, are depicted in Figure \ref{Chords44}.

We are going to find the distribution of the length of the chords of the sets of the form $D(C) \cap Cr(a)$. For each direction $\theta$ of the chords, the density of measure of $D(C, \theta) \cap Cr(a, \theta)$ and $D(C, \theta) \cap Cr(b, \theta)$ is a step function of the same width as the density of measure of $D(C, \theta)$, but with shorter heights which add up to the height of the density of measure of $D(C, \theta)$. In Figure \ref{Chords1} we saw that this height was $c \sin\theta$ (formula (\ref{1})). In the same figure we can see that, in an obvious notation,

\begin{equation}
\mu(D(C, \theta) \cap Cr(b, \theta)) = b \sin(\theta-\alpha)
\end{equation}
and
\begin{equation}
\mu(D(C, \theta) \cap Cr(a, \theta)) = a \sin(\pi-\beta-\theta).
\end{equation}

We now follow the steps of section 1 changing $\mu_C(\theta)$ and $\mu_A(\theta)$ by the above expressions.

For $D(C, \theta) \cap Cr(b, \theta)$ in formula (5) we substitute $\mu_C(\theta)$ by $\mu(D(C, \theta) \cap Cr(b, \theta))$, and in formulae (9) and (10) we substitute their right hand sides by ${b c \over 2 A} \sin\theta \sin(\theta-\alpha)$ and ${1 \over 2} \Big( (bb-aa) \cos \alpha + \sin(aa-bb) \cos (aa+bb-\alpha) \Big)$, respectively. We prefer to call the two segments which meet $a$ and $b$, rather than $b$ and $c$. Therefore we do the substitution $(a, \alpha) \rightarrow (c, \gamma), (b, \beta) \rightarrow (a, \alpha)$ and $(c, \gamma) \rightarrow (b, \beta)$ to find that the density of measure in $Cr(a) \cap D(B)$ is

$$
\rho_{aB}(\ell) = {a b \over 4 A} \left[ H\left({2 A \over b} - \ell \right) \big( \beta \cos\gamma  + \cos\alpha\ \sin\beta \big) \right. +
$$

$$
H(c - \ell) H\bigg(\ell - {2 A \over b}\bigg) \Bigg( \bigg(\arcsin{2 A \over b \ell} - {2 A \over b \ell} \sqrt{1 - \left({2 A \over b \ell}\right)^2} - \alpha \bigg) \cos\gamma + \bigg({2 A \over b \ell}\bigg)^2 \sin\gamma - \sin\alpha \cos\beta \bigg) +
$$

\begin{equation}
\left.
H(a - \ell) H\bigg(\ell - {2 A \over b}\bigg) \Bigg( \bigg(\arcsin{2 A \over b \ell} - {2 A \over b \ell} \sqrt{1 - \left({2 A \over b \ell}\right)^2} - \gamma \bigg) \cos\gamma + \bigg(1 - \bigg({2 A \over b \ell}\bigg)^2 \bigg) \sin\gamma \Bigg) \right].
\end{equation}

When the height $h_B$ does not hit the side $b$, as shown in Figure \ref{Chords555}, there are only two terms:

\begin{figure}[ht]
  \begin{center}
  \includegraphics[width=0.4 \textwidth]{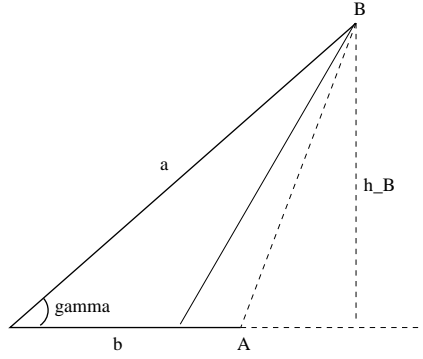}\caption{\label{Chords555}The height $h_B$ does not hit the side $b$}
  \end{center}
\end{figure}

$$
\rho_{aB}(\ell) = {a b \over 4 A} \left[ H\left(c - \ell \right) \big( \beta \cos\gamma  + \cos\alpha\ \sin\beta \big) \right. +
$$

$$
H(\ell - c) H(a - \ell) \Bigg( \bigg(\arcsin{2 A \over b \ell} - {2 A \over b \ell} \sqrt{1 - \left({2 A \over b \ell}\right)^2} - \gamma \bigg) \cos\gamma + \bigg( 1 - \bigg({2 A \over b \ell}\bigg)^2 \bigg) \sin\gamma \bigg).
$$

\begin{figure}[ht]
  \begin{center}
  \includegraphics[width=0.65 \textwidth]{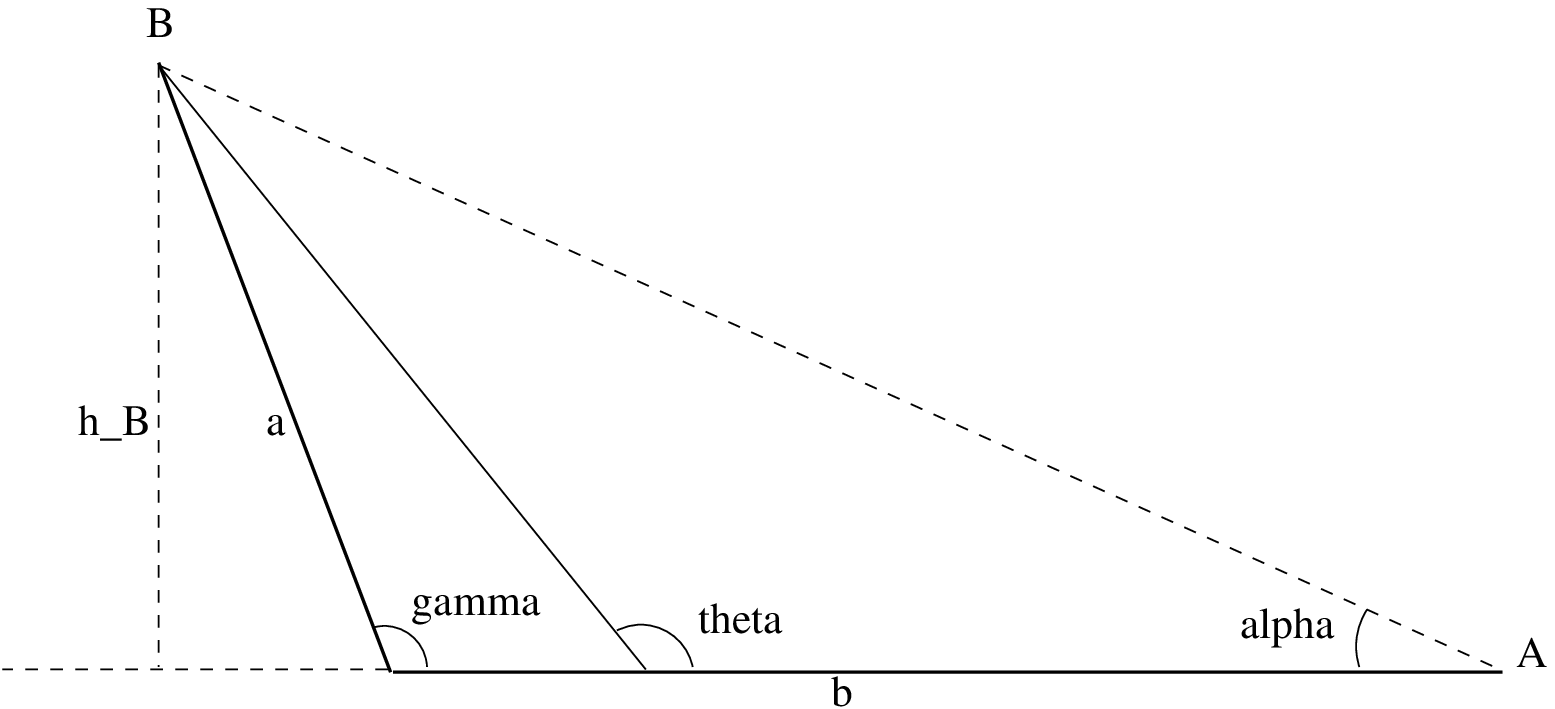}\caption{\label{Chords6666}$\gamma > {\pi \over 2}$}
  \end{center}
\end{figure}

When $\gamma > {\pi \over 2}$ (Figure \ref{Chords6666}), a similar calculation yields

$$
\rho_{aB}(\ell) =
$$

$$
{a b \over 2 A} H(a - \ell) \int_{\gamma}^{\pi - \alpha} d\theta\ \sin\theta \sin(\theta-\gamma) +
{a b \over 2 A} H(\ell - a) H(c - \ell) \int_{\pi - \arcsin{h_A \over \ell}}^{\pi - \alpha} d\theta\ \sin\theta \sin(\theta-\gamma) =
$$

$$
{a b \over 4 A} \left[ H\left(a - \ell \right) \big( \beta \cos\gamma  + \cos\alpha\ \sin\beta \big) \right. +
$$

$$
H(\ell - a) H(c - \ell) \Bigg( \bigg(\arcsin{2 A \over b \ell} - {2 A \over b \ell} \sqrt{1 - \left({2 A \over b \ell}\right)^2} + \sin \alpha \cos \alpha- \alpha \bigg) \cos\gamma + \bigg( \bigg({2 A \over b \ell}\bigg)^2 - (\sin \alpha)^2 \bigg) \sin\gamma \bigg).
$$

Since

\begin{equation}
Cr(a) \cap Cr(b) = (Cr(a) \cap D(B)) \sqcup (Cr(b) \cap D(A)),
\end{equation}
to find the density of measure of the length of the chords which cross the sides $a$ and $b$ we have to add the corresponding densities of the chords in the sets $Cr(a) \cap D(B)$ and $Cr(b) \cap D(A)$.

We have proven

{\bf{Proposition 3.1.}} {\it{The density of measure of the length of chords which cut two concurring segments of lengths $a$ and $b$ and which form an angle $\gamma$ is given by}}

\begin{equation}
\rho(a, b, \gamma; \ell) \equiv
\rho_{aB}(\ell, \alpha, \beta, \gamma, a, b, c) + \rho_{bA}(\ell, \alpha, \beta, \gamma, a, b, c).
\end{equation}

\section{Distribution of the length of chords between two segments}

{\bf{Proposition 4.1.}} {\it{The density of measure of the length of chords which cut two non parallel segments of
extremes $A, A'$ and $B, B'$ and which form an angle $\gamma \neq 0$ is given at $\ell$ by

$$
\rho(A, A', B, B'; \ell) =
$$
\begin{equation}
\rho(OA', OB', \gamma; \ell) - \rho(OA, OB', \gamma; \ell) - \rho(OA', OB, \gamma; \ell) + \rho(OA, OB, \gamma; \ell).
\end{equation}
where the $\rho$'s in the rhs denote the density of measure of Proposition 3.1 and OA, OA', OB, OB' denote lengths.}}

\begin{figure}[ht]
  \begin{center}
  \includegraphics[width=0.45 \textwidth]{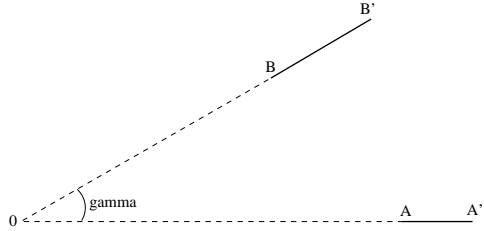}\caption{\label{Chords55}Any two non parallel segments define an angle}
  \end{center}
\end{figure}

{\bf{Proof.}} It is clear from Figure \ref{Chords55} that
$$
\rho(OA', OB', \gamma; \ell) = \rho(OA, OB', \gamma; \ell) + \rho(OA', OB, \gamma; \ell) - \rho(OA, OB, \gamma; \ell) + \rho(A, A', B, B'; \ell),
$$
from which the result follows immediately. $\Box$
\medskip

The parallel segments case can be found as a limit of the previous case, but it is faster to obtain it directly, because when the segments are parallel there are, for a given length, at most two directions for which the chords have that given length.

\begin{figure}[ht]
  \begin{center}
  \includegraphics[width=0.7 \textwidth]{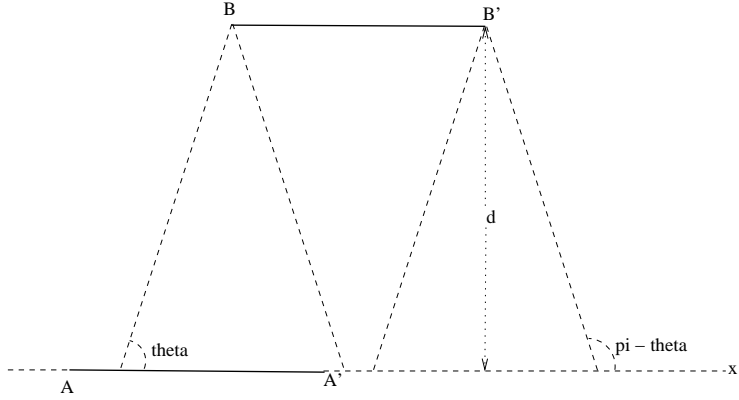}\caption{\label{Chords666}Two parallel segments.}
  \end{center}
\end{figure}

In Figure \ref{Chords666} we see two segments whose supporting lines are a distance $d$ apart. The density of measure of the chords of length $\ell = d$ is obviously proportional to the length of the intersection between $AA'$ and the projection of $BB'$ onto the $x$ axis, that is,

\begin{equation}
\rho(d) \propto h\big( \min(x_{B'}, x_{A'}) - \max(x_{B}, x_{A}) \big),
\end{equation}
where
\begin{equation}
h(x) \equiv x H(x).
\end{equation}
When $\ell > d$ there are two directions for which the chords have length $\ell$. These directions are $\arcsin {d \over \ell}$ and $\pi - \arcsin {d \over \ell}$. The density of measure of the chords of length $\ell$ and direction $\arcsin {d \over \ell}$ (denoted $\theta$ in Figure \ref{Chords6666}) is proportional to the length of the intersection between $AA'$ and the projection of $BB'$ along the lines of direction $\arcsin {d \over \ell}$ times ${d \over \ell}$. A factor ${d \over \ell}$ is needed to convert the length of the intersection into the cross section of the set of chords. Finally, an angle dependent factor is needed to compute the density of measure. Since

\begin{equation}
\theta = \arcsin {d \over \ell}\ \ \Rightarrow\ \ d\theta = -{d \over \ell \sqrt{\ell^2 - d^2}} d\ell,
\end{equation}
this factor is ${d \over \ell \sqrt{\ell^2 - d^2}}$. The measure of the chords of length $\ell$ and direction $\pi - \arcsin {d \over \ell}$ can be computed similarly.

$$
\rho(\ell) = {d^2 \over \ell^2 \sqrt{\ell^2 - d^2}} \bigg[ h\big( \min(x_{B'} - \sqrt{\ell^2 - d^2}, x_{A'}) - \max(x_{B} - \sqrt{\ell^2 - d^2}, x_{A}) \big) +
$$

\begin{equation}
h\big( \min(x_{B'} + \sqrt{\ell^2 - d^2}, x_{A'}) - \max(x_{B} + \sqrt{\ell^2 - d^2}, x_{A}) \big)  \bigg],
\ \ \ \ell > d.
\end{equation}
If the previous definition is taken to hold also when $\ell = d$, only a measure zero error is made.

We have proven

{\bf{Proposition 4.2.}} {\it{The density of measure of the length of chords which cut two parallel segments whose supporting lines are at a distance $d$ is given by the previous expression when $\ell \geq d$ and is zero otherwise.}}

\bigskip

\section{Distribution of chords in a convex polygon}

In order to compute the first formula of the preceding section when $A, A', B$ and $B'$ are given in Cartesian coordinates we need to compute, first of all, the coordinates of $O$. It is tedious but elementary to find that the straight lines supporting the segments $AA'$ and $BB'$ meet at a point $O$ of coordinates

\begin{equation}
(x_O, y_O) =
\left(
{
\left|
  \begin{array}{cc}
    \left|
      \begin{array}{cc}
      x_{A}  & y_{A}  \\
      x_{A'} & y_{A'} \\
      \end{array}
    \right| &  x_{A'} - x_{A} \\
     &  \\
    \left|
      \begin{array}{cc}
      x_{B}  & y_{B}  \\
      x_{B'} & y_{B'} \\
      \end{array}
    \right|
     &  x_{B'} - x_{B} \\
  \end{array}
\right|
\over
\left|
  \begin{array}{cc}
    x_{A'} - x_{A} &  y_{A'} - y_{A} \\
    x_{B'} - x_{B} &  y_{B'} - y_{B} \\
  \end{array}
\right|
}
,
{
\left|
  \begin{array}{cc}
    \left|
      \begin{array}{cc}
      x_{A}  & y_{A}  \\
      x_{A'} & y_{A'} \\
      \end{array}
    \right| &  y_{A'} - y_{A} \\
     &  \\
    \left|
      \begin{array}{cc}
      x_{B}  & y_{B}  \\
      x_{B'} & y_{B'} \\
      \end{array}
    \right|
     &  y_{B'} - y_{B} \\
  \end{array}
\right|
\over
\left|
  \begin{array}{cc}
    x_{A'} - x_{A} &  y_{A'} - y_{A} \\
    x_{B'} - x_{B} &  y_{B'} - y_{B} \\
  \end{array}
\right|
}
\right)
\end{equation}

The lengths which appear in formulae (22) (to which formula (23) refers) and (23) are $OA, OA', OB, OB'$ and $A'B'$, all of which are straightforward to compute because the coordinates of the five points involved are now known. As for the angles of the triangle $OA'B'$, they are

\begin{equation}
\gamma =
\arccos {\overrightarrow{AA'} \cdot \overrightarrow{BB'} \over |\overrightarrow{AA'}| |\overrightarrow{BB'}|},
\end{equation}
and, by the cosine theorem,
\setcounter{equation}{31}
$$
\alpha =
\arccos {|\overrightarrow{BB'}| - |\overrightarrow{AA'}| \cos \gamma \over \sqrt{|\overrightarrow{AA'}|^2 + |\overrightarrow{BB'}|^2 - 2 |\overrightarrow{AA'}| |\overrightarrow{BB'}| \cos \gamma}}
$$
and

\begin{equation}
\beta =
\arccos {|\overrightarrow{AA'}| - |\overrightarrow{BB'}| \cos \gamma \over \sqrt{|\overrightarrow{AA'}|^2 + |\overrightarrow{BB'}|^2 - 2 |\overrightarrow{AA'}| |\overrightarrow{BB'}| \cos \gamma}}.
\end{equation}
The ``missing" side of the triangle can also be found using the cosine theorem. It is now possible to write $\rho(A, A', B, B'; \ell)$ in Cartesian coordinates.

{\bf{Proposition 5.1.}} {\it{The density of measure of the length of chords of a convex polygon given by the Cartesian coordinates of its vertices can be found using the formulae of sections 4 and 5.}}

{\bf{Proof.}} In a convex polygon, up to measure zero, all chords go trough exactly two sides. If a convex polygon is given by the Cartesian coordinates of its vertices, then its total measure density is the sum of the densities of all pairs of sides, which the formulae of this and the preceding section show how to find. $\Box$

That the set of chords cutting a convex polygon can be decomposed in the (up to measure zero) disjoint union of the sets of chords cutting the pairs of sides is a straightforward proposition which has been known and used by other authors (e. g. \cite{Mallows1970, Harutyunyan2011}). The advantage of the formulae presented up to here over the formulae presented by H. S. Harutyunyan and V. K. Ohanyanian \cite{Harutyunyan2011} is that the formulae of sections 4 and 5 allow for a straightforward substitution of the Cartesian coordinates (cf. \cite{Harutyunyan2011}, section 5).
The author has written a file of Mathematica code that does just that and which can be downloaded at \cite{RGPN}.

\section{Convex examples}

\begin{figure}[ht]
  \begin{center}
  \includegraphics[width=0.75 \textwidth]{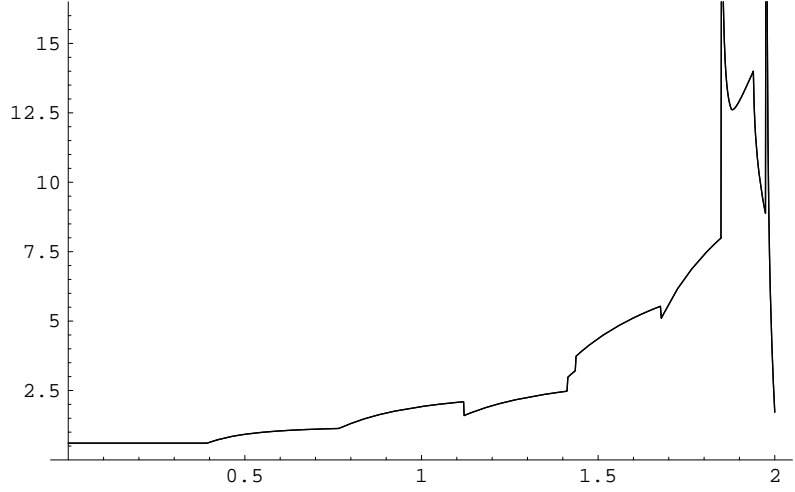}\caption{\label{Dodecagono}Measure density of a pair of Mallows and Clark dodecagons}
  \end{center}
\end{figure}

We have shown how to compute the distribution of chords of a (convex) polygon by computing it for each pair and then adding all the distributions found in this manner. This leads naturally to the question of whether a couple of non-congruent $n$-gons exists such that the $n \choose 2$ pairs of sides are, setwise, congruent to each other. Such a pair (actually, a whole family of pairs) was found by Mallows and Clark \cite{Mallows1970} in 1970, and we show in Figure \ref{Dodecagono} their measure density.

The particular pair taken as example is based on the octagon of coordinates $\{A, B, C, D,\\ E, F, G, H \} \equiv \{ (1, 0),\ (1/\sqrt{2},1/\sqrt{2}),\ (0,1),\ (-1/\sqrt{2},1/\sqrt{2}),\ (-1,0),\ (-1/\sqrt{2},-1/\sqrt{2}),\\ (0,-1),\ (1/\sqrt{2},-1/\sqrt{2}) \}$. To obtain one of the dodecagons we add the vertices $\{ AB, BC,\\ DE, EF \} \equiv \{ 1.1 {A+B \over 2}, 1.1 {B+C \over 2}, 1.1 {D+E \over 2}, 1.1 {E+F \over 2} \}$, and to obtain the other the vertices $\{ AB, DE, EF, FG \} \equiv \{ 1.1 {A+B \over 2}, 1.1 {D+E \over 2}, 1.1 {E+F \over 2}, 1.1 {F+G \over 2} \}$. Although not clear in Figure \ref{Dodecagono}, the graph is continuous, since the formulae obtained in the preceding section are continuous. More generally, this follows from a theorem of Sulanke for bounded, non-degenerate, convex figures (\cite{Sulanke1961}, page 55, Proposition 2). This particular example is worked out in the file Convex.nb in \cite{RGPN}. With that code the above graph is produced in about 30 seconds on a personal computer.

\section{The concave quadrilateral}

\begin{figure}[ht]
  \begin{center}
  \includegraphics[width=0.6 \textwidth]{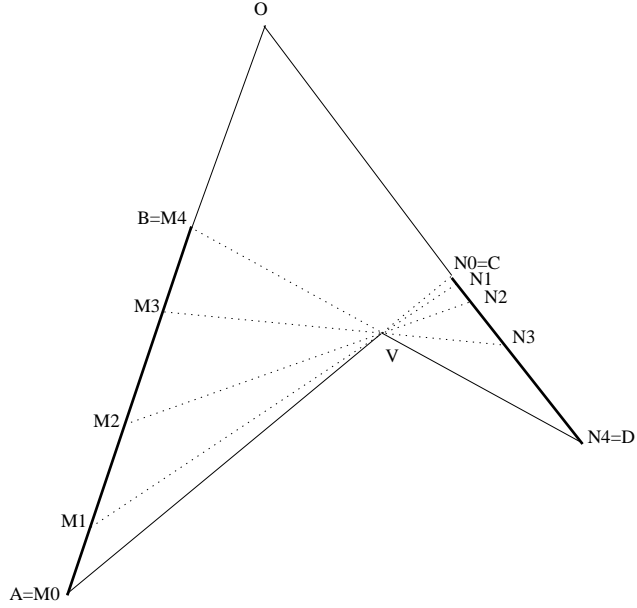}\caption{\label{PrePajarita}The concave quadrilateral.}
  \end{center}
\end{figure}

The measure density for the depicted quadrilateral can be decomposed in measure densities of pairs of segments as follows:

\begin{equation}
\rho_{(V A, A O)} + \rho_{(V D, D O)} + \rho_{(V A, O C)} + \rho_{(V D, O B)} + \rho_{(A O, O C)} + \rho_{(B O, O D)} - \rho_{(O B, O C)} + \rho_{conc(A B, C D)}.
\end{equation}
Of these eight terms the first seven are of the kind already investigated, while the last one is new, because not every segment joining $AB$ and $CD$ is contained in the quadrilateral.

The new term can be approximated (and bounded by above) by
$$
\rho(M_0, M_1, N_0, N_1; \ell) + \rho(M_1, M_2, N_0, N_2; \ell) + \rho(M_2, M_3, N_0, N_3; \ell) + \rho(M_3, M_4, N_0, N_4; \ell) =
$$

$$
\rho(M_0, M_0 + m_1 \vec u, N_0, op\left( M_0 + m_1 \vec u \right); \ell) + \rho(M_0 + m_1 \vec u, M_0 + m_2 \vec u, N_0, op\left( M_0 + m_2 \vec u \right); \ell)
$$

\begin{equation}
+ \rho(M_0 + m_2 \vec u, M_0 + m_3 \vec u, N_0, op\left( M_0 + m_3 \vec u \right); \ell) + \rho(M_0 + m_3 \vec u, M_0 + m_4 \vec u, N_0, op\left( M_0 + m_4 \vec u \right); \ell),
\end{equation}

\noindent
where $\vec u = {B - A \over |B - A|},\ m_i = |M_i - M_0|$, $op$ is a geometrical function which assigns to a point in the left side the point in the right side such that the straight line which they determine contains the vertex $V$ and $\rho$ is the measure density function for two segments which has already been studied. The function $op$ is derived in the Appendix.

For any $n \in N$ there is a set of equally spaced points $\{ M_0, M_1, \cdots, M_n \}$ in the right side and an approximation as above can be written:

\begin{equation}\label{35}
\sum_{i=1}^n \rho(M_0 + {i-1 \over n} |M_1 - M_0| \vec u, M_0 + {i \over n} |M_1 - M_0| \vec u, N_0, op\left( M_0 + {i \over n} |M_1 - M_0| \vec u \right); \ell).
\end{equation}

This expression is the starting point for two different ways of computing the measure density for a concave quadrilateral. The first way to be presented is numerically more efficient, while the second one is mathematically more sophisticated and allows, sometimes, an analytical computation of the measure density.

In section 10 we shall refer repeatedly to a pentagon of the kind $ABCDV$. For this reason we give the following definition.

{\bf{Definition 7.1}} {\it{A concave pentagon whose concave vertex is the intersection of the diagonals of the convex quadrilateral whose vertices are the other four points will be called a diagonal pentagon.}}

Note that there are four diagonal pentagons for a given set of five vertices satisfying the above definition. Note also that in a diagonal pentagon the distribution of the length of chords between any pair of sides can be computed using Propositions 4.1 and 4.2, except for one pair, which is the pair $AB, CD$ in the case of Figure \ref{PrePajarita}.

\section{Converging bounds for the measure density of the concave quadrilateral}\label{Converging bounds}

As remarked in the last section, the last expression is an upper bound for the measure density $\rho_{conc(A B, C D)}$. This is better seen in the case $n=4$, depicted in Figure \ref{PrePajarita}. In the expression (\ref{35}) all chords joining the segments $AB$ and $CD$ and passing above the vertex are counted in, as well as some chords which pass below the vertex $V$. Similarly, the expression

\begin{equation}\label{36}
\sum_{i=1}^n \rho(M_0 + {i-1 \over n} |M_1 - M_0| \vec u, M_0 + {i \over n} |M_1 - M_0| \vec u, N_0, op\left( M_0 + {i-1 \over n} |M_1 - M_0| \vec u \right); \ell).
\end{equation}
is, for any $n \in N$, a lower bound for the measure density $\rho_{conc(A B, C D)}$.

Using the notation $app_{-}(n)(\ell)$ and $app_{+}(n)(\ell)$ for the last two sums, we have remarked that

\begin{equation}
app_{-}(n)(\ell) \leq \rho_{conc}(AB, CD)(\ell) \leq app_{+}(n)(\ell),\ \ \ \ \forall n \in N.
\end{equation}
The difference of the sums is

$$
app_{+}(n)(\ell) - app_{-}(n)(\ell) =
$$

$$
\sum_{i=1}^n \rho(M_0 + {i-1 \over n} |M_1 - M_0| \vec u, M_0 + {i \over n} |M_1 - M_0| \vec u,
$$

\begin{equation}
op\left( M_0 + {i-1 \over n} |M_1 - M_0| \vec u \right), op\left( M_0 + {i \over n} |M_1 - M_0| \vec u \right); \ell)
\ \underset{n \rightarrow \infty}\longrightarrow\ 0.
\end{equation}
This is intuitive because the difference $app_{+}(n)(\ell) - app_{-}(n)(\ell)$ involves, as $n \rightarrow \infty$, only the chords that go through the vertex $V$, which is a set of measure zero. It can be checked rigorously using the analytic form of the density of measure for the chords between two segments.

{\bf{Proposition 8.1.}} {\it{The last two equations prove that}}

\begin{equation}
\lim_{n \rightarrow \infty} app_{-}(n)(\ell) = \lim_{n \rightarrow \infty} app_{+}(n)(\ell)
= \rho_{conc}(AB, CD)(\ell).
\end{equation}

\section{The measure density of the concave quadrilateral as a Riemann integral}
We momentarily omit the arguments $N_0$ and $\ell$ of $\rho$ in the sum (35) and use the following short hand notation for the other three:

\begin{equation}
\sum_{i=1}^n \rho(i-1, i, op\left( i \right)).
\end{equation}
$\rho$ is $C^2$ except at a finite number of points (see formulae in the sections on convex polygons), therefore:

$$
\sum_{i=1}^n \rho(i-1, i, op\left( i \right)) =
$$

\begin{equation}
\sum_{i=1}^n \rho(i-1, i-1, op\left( i \right)) + \rho'_2(i-1, i-1, op\left( i \right))\ {|B - A| \over n} + o\left( {1 \over n} \right)
\end{equation}
where $\rho'_2$ is the partial derivative with respect to the second argument. The term $\rho(i-1, i-1, op\left( i \right))$ is zero because the segment of extremes $M_0 + {i-1 \over n} |B - A| \vec u$ and $M_0 + {i-1 \over n} |B - A| \vec u$ is of zero length. The expression can be further transformed as follows:

$$
\sum_{i=1}^n \rho'_2(i-1, i-1, op\left( i \right))\ {|B - A| \over n} + o\left( {1 \over n} \right) =
$$

$$
\sum_{i=1}^n \rho'_2(i-1, i-1, op\left( i-1 \right))\ {|B - A| \over n} + o\left( {1 \over n} \right) +
$$

$$
\rho''_{23}(i-1, i-1, op\left( i-1 \right)) op'(i-1)\ {|B - A|^2 \over n^2} + o\left( {1 \over n^2} \right) =
$$

\begin{equation}
\sum_{i=1}^n \rho'_2(i-1, i-1, op\left( i-1 \right))\ {|B - A| \over n} + o\left( {1 \over n} \right).
\end{equation}

Had we started from expression (\ref{36}) we would have also reached expression (42). Since expression (\ref{36}) is, for any $n$, always smaller than or equal to expression (\ref{36}) and the differ by at most $\sum_{i=1}^n o\left( {1 \over n} \right)$, which satisfies

\begin{equation}
\lim_{n \to \infty} \sum_{i=1}^n o\left( {1 \over n} \right) = \lim_{n \to \infty} n\ o\left( {1 \over n} \right) =
\lim_{n \to \infty} n\ {1 \over n} {o\left( {1 \over n} \right) \over {1 \over n}} =
\lim_{n \to \infty} {o\left( {1 \over n} \right) \over {1 \over n}} = 0,
\end{equation}
then the limit of expression (42) as $n \rightarrow \infty$ is the following Riemann integral:

\begin{equation}
\int_{0}^{|B - A|} da\ \rho'_2\left(a, a, op(a) \right) =
\int_{0}^{|B - A|} da\  \rho'_2(A + a \vec u, A + a \vec u, C, op\left( A + a \vec u \right); \ell),
\end{equation}
where in the rhs we drop the short hand notation and $A$, $B$ and $C$ are the extremes of the thick intervals depicted in Figure \ref{PrePajarita}.

We have proven

{\bf{Proposition 9.1.}}
\begin{equation}\label{45}
\int_{0}^{|B - A|} da\  \rho'_2(A + a \vec u, A + a \vec u, C, op\left( A + a \vec u \right); \ell)
= \rho_{conc}(AB, CD)(\ell).
\end{equation}

In the remaining part of the section a formula for the integrand of the above expression will be derived. Note that the approximations $app_{-}(n)(\ell)$ and $app_{+}(n)(\ell)$ discussed in the previous section are lower and upper bounds of $\rho'_2(\ell)$, but not lower and upper Riemann integrals of the preceding Riemann integral.

\bigskip
\subsection{The derivative $\rho'_2$}

We need to find the derivative $\rho'_2$. In this subsection we use the notation of section 3. First we address the general case in which the two segments are not parallel. From (23) it is clear that

\begin{equation}
\rho'_2(A, A', B, B'; \ell) = \rho'_1(OA', OB', \gamma; \ell) - \rho'_1(OA', OB, \gamma; \ell).
\end{equation}
Since
\setcounter{equation}{21}

\begin{equation}
\rho(a, b, \gamma; \ell) =
\rho_{aB}(\ell, \alpha, \beta, \gamma, a, b, c) + \rho_{bA}(\ell, \alpha, \beta, \gamma, a, b, c),
\end{equation}
\setcounter{equation}{46}
we need to find the {\it{total}} derivatives ${d \over da} \rho_{aB}(\ell, \alpha, \beta, \gamma, a, b, c)$ and ${d \over da} \rho_{bA}(\ell, \alpha, \beta, \gamma, a, b, c)$. The first of these total derivatives is

$$
{d \over da} \rho_{aB}(\ell, \alpha, \beta, \gamma, a, b, c) = \rho'_{aB2}(\ell, \alpha, \beta, \gamma, a, b, c) {\partial \alpha \over \partial a} da + \rho'_{aB3}(\ell, \alpha, \beta, \gamma, a, b, c) {\partial \beta \over \partial a} da +
$$

\begin{equation}
\rho'_{aB5}(\ell, \alpha, \beta, \gamma, a, b, c) da + \rho'_{aB7}(\ell, \alpha, \beta, \gamma, a, b, c) {\partial c \over \partial a} da.
\end{equation}
From the cosine theorem, $c = \sqrt{a^2 + b^2 - 2 a b \cos \gamma}$,

\begin{equation}
{\partial c \over \partial a} = {a - b \cos \gamma \over c}.
\end{equation}
From the sine theorem, ${\sin \alpha \over a} = {\sin \gamma \over \sqrt{a^2 + b^2 - 2 a b \cos \gamma}} \Rightarrow \alpha = \arcsin a {\sin \gamma \over \sqrt{a^2 + b^2 - 2 a b \cos \gamma}} $,

\begin{equation}
{\partial \alpha \over \partial a} = {b \sin \gamma \over c^2}.
\end{equation}
Likewise,

\begin{equation}
{\partial \beta \over \partial a} = -{b \sin \gamma \over c^2}.
\end{equation}

A similar calculation is required for ${d \over da} \rho_{bA}(\ell, \alpha, \beta, \gamma, a, b, c)$. This approach leads to an expression for $\rho'_2$ which is a few pages long and, therefore, of questionable advantage with respect to simulations.
\bigskip

\begin{figure}[ht]
  \begin{center}
  \includegraphics[width=0.6 \textwidth]{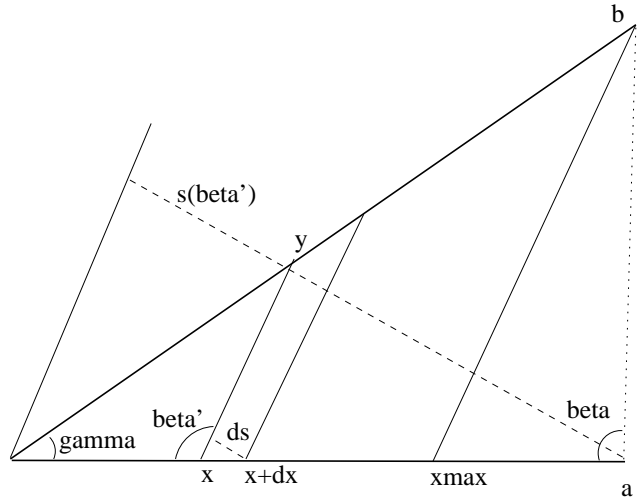}\caption{\label{Conc2}In the example above $\beta \approx \pi/2.\ x_{max}(\beta')$ is, for a given $\beta'$, the maximum $x$ for which the chord still cuts the upper side.}
  \end{center}
\end{figure}

\subsection{A more direct calculation of $\rho'_2$. Preliminaries}

$\rho'_2$ can be found by a more direct method, which is to consider the distribution of the length of chords between two segments, one of which is of infinitesimal length.

As a first step towards this goal we are going to find the measure of chords that go through two concurring segments as a double integral over the segments. Let $x$ and $y$ be coordinates along the said segments, as shown in Figure \ref{Conc2}. To find the said measure we find the section of the chords for a given direction and then integrate over the directions:

$$
\mu = \int_0^{\pi-\gamma} d\beta'\ s(\beta') = \int_0^{\pi-\gamma} d\beta'\ \int_0^{s(\beta')} ds =
$$

$$
\int_0^{\beta} d\beta' \int_0^a dx\ {\partial s(\beta', x) \over \partial x} + \int_{\beta}^{\pi-\gamma} d\beta' \int_0^{x_{max}(\beta')} dx\ {\partial s(\beta', x) \over \partial x} =
$$

\begin{equation}
\int_0^{\beta} d\beta' \int_0^a dx\ \sin \beta' + \int_{\beta}^{\pi-\gamma} d\beta' \int_0^{x_{max}(\beta')} dx\ \sin \beta',
\end{equation}
where $\beta$ is the angle between the horizontal side and the missing side of the triangle and

\begin{equation}
ds = \sin \beta' dx,
\end{equation}
which is obvious from Figure \ref{Conc2}, has been used.

\begin{equation}
\mu = \int_0^a dx\ \int_0^{\beta} d\beta' \sin \beta' + \int_0^{a} dx\ \int_{\beta}^{\beta'_{max}(x)} d\beta' \sin \beta' = \int_0^{a} dx\ \int_{0}^{\beta'_{max}(x)} d\beta' \sin \beta'.
\end{equation}

\begin{figure}[ht]
  \begin{center}
  \includegraphics[width=0.4 \textwidth]{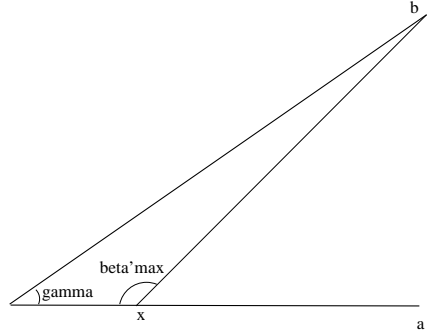}\caption{\label{Conc3}For a given $x$, there is a maximum $\beta'_{max}(x)$ for which the chord still cuts the upper side.}
  \end{center}
\end{figure}

To change the integration variable in the last integral from $\beta'$ to $y$ we use the sine theorem in the form

$$
{\sin \beta' \over y} = {\sin \gamma \over \sqrt{x^2 + y^2 - 2 x y \cos \gamma}}\ \ \Rightarrow\ \
\beta' = \arcsin {y \sin \gamma \over \sqrt{x^2 + y^2 - 2 x y \cos \gamma}}\ \ \Rightarrow\ \
$$

\begin{equation}
{d \beta' \over dy} = {x \sin \gamma \over x^2 + y^2 - 2 x y \cos \gamma}.
\end{equation}
Substitution in the last integral yields

\begin{equation}
\mu = \int_0^{a} \int_{0}^{b} dx\ dy\ {x y (\sin \gamma)^2 \over (x^2 + y^2 - 2 x y \cos \gamma)^{3/2}}.
\end{equation}
The above formula holds for two concurring segments.

To take the $\gamma \rightarrow 0$ limit in the above integrand we note that the denominator is the cube of the distance between the points of coordinates $x$ and $y$ and that the numerator is the product of the heights of the triangle which go through the points of coordinates $x$ and $y$. Then, when the lines are parallel, the above integrand becomes

\begin{equation}\label{56}
{d^2 \over (d^2 + (y-x)^2)^{3/2}},
\end{equation}
where $d$ is the distance between the parallel lines.

\subsection{A more direct calculation of $\rho'_2$. Conclusion}

It follows from the equality (\ref{45}) that $\rho'_2(A + a \vec u, A + a \vec u, C, op\left( A + a \vec u \right); \ell)\ da$ is the measure density of $\ell$ for the chords which join the segments $A + a \vec u, A + (a + da) \vec u$ and $C, op\left( A + a \vec u \right)$. Therefore, using a notation that keeps only the last argument of $\rho'_2$, we may compute it as follows:

\begin{figure}[ht]
  \begin{center}
  \includegraphics[width=0.4 \textwidth]{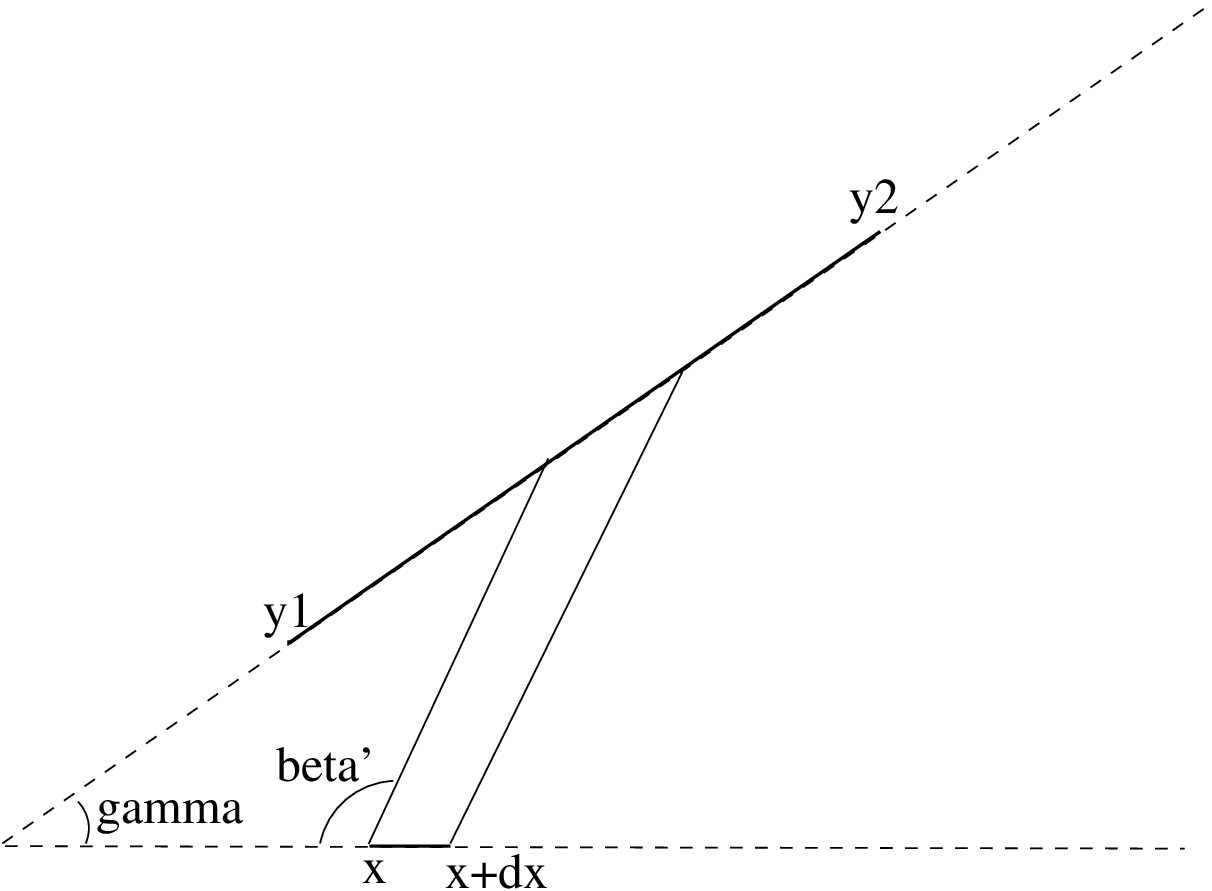}\caption{\label{Conc4}}
  \end{center}
\end{figure}

\begin{equation}
\rho'_2(\ell) = \lim_{dx \rightarrow 0} {1 \over dx} \int_{y_1}^{y_2}dy\ \int_{x}^{x+dx}dx\ {x y (\sin \gamma)^2 \over (x^2 + y^2 - 2 x y \cos \gamma)^{3/2}}\ \delta(\sqrt{x^2 + y^2 - 2 x y \cos \gamma} - \ell),
\end{equation}
where the coordinates $x$ and $y$ have their origin in the vertex shown in Figure \ref{Conc4}. The comparison of Figures 10 and 13 yields the relations between the arguments of $\rho'_2$ in (\ref{45}) and the arguments of $\rho'_2$ used above:

\begin{equation}\label{58}
|OC| \rightarrow y_1,\ \ \ |O op\left( A + a \vec u \right)| \rightarrow y_2,\ \ \ |OA - a \vec u| \rightarrow x.
\end{equation}

To take advantage of the Dirac delta function we have to change the variables of integration from $(x,y)$ to $\phi(x,y) \equiv \sqrt{x^2 + y^2 - 2 x y \cos \gamma}$ and some other variable. We choose the other one to be the angle $\beta'$. By the sine theorem, the relation between the variables is:

\begin{equation}
\begin{array}{ccc}
  x & = & {\phi \over \sin \gamma} \sin(\gamma + \beta') \\
  & &  \\
  y & = & {\phi \over \sin \gamma} \sin \beta'\ \ \ \ \ \ \ \\
\end{array}
\end{equation}
The integrand becomes

$$
\left|{\partial(x,y) \over \partial(\phi, \beta')} \right| {x y (\sin \gamma)^2 \over (x^2 + y^2 - 2 x y \cos \gamma)^{3/2}}\ \delta(\sqrt{x^2 + y^2 - 2 x y \cos \gamma} - \ell) =
$$

\begin{equation}
{\phi \over \sin\gamma} {\sin\beta' \sin(\beta'+\gamma) \over \phi} \delta(\phi - \ell) =
{\sin\beta' \sin(\beta'+\gamma) \over \sin\gamma} \delta(\phi - \ell).
\end{equation}
We are going to integrate first over $\phi$. For fixed $\beta'$ the limits of integration of $\phi$ are, by the sine theorem, ${\sin\gamma \over \sin(\beta'+\gamma)}\ x$ and ${\sin\gamma \over \sin(\beta'+\gamma)}\ (x + dx)$. Therefore, in the expression (where the integration limits have been omitted)

\begin{equation}
\rho'_2(\ell) = \lim_{dx \rightarrow 0} {1 \over dx} \int d\beta'\ {\sin\beta' \sin(\beta'+\gamma) \over \sin\gamma} \int d\phi\ \delta(\phi - \ell),
\end{equation}
the integral over $\phi$ will be 1 when $\ell \in \left[ {\sin\gamma \over \sin(\beta'+\gamma)}\ x,\ {\sin\gamma \over \sin(\beta'+\gamma)}\ (x + dx) \right]$ and 0 otherwise. The way to capture this into the above expression is to give appropriate limits of integration in the integral over $\beta'$. To do this we note that

\begin{equation}
\ell \in \left[ {\sin\gamma \over \sin(\beta'+\gamma)}\ x,\ {\sin\gamma \over \sin(\beta'+\gamma)}\ (x + dx) \right]
\end{equation}
is equivalent to

$$
\beta' \in \left[\arcsin\left( {x \over \ell} \sin\gamma \right) - \gamma,\ \arcsin\left( {x + dx \over \ell} \sin\gamma \right) - \gamma \right] \cup
$$

\begin{equation}
\left[\pi - \arcsin\left( {x + dx \over \ell} \sin\gamma \right) - \gamma,\ \pi - \arcsin\left( {x \over \ell} \sin\gamma \right) - \gamma \right] = I_1 \cup I_2,
\end{equation}
where $I_1$ and $I_2$ are notations for the intervals above and by $\arcsin$ the branch that takes values between $-\pi/2$ and $\pi/2$ is meant. The two intervals correspond to the two chords that have a given length $\ell$, one to each side of the shortest segment that goes from the point of coordinate $x$ to the straight line where the point of coordinate $y$ lies. Let the angles $\beta'$ which correspond to the coordinates $y_1$ and $y_2$ be $\beta_1$ and $\beta_2$ (which satisfy $\sin \beta_1 = {y_1 \over \ell} \sin \gamma$ and $\sin \beta_2 = {y_2 \over \ell} \sin \gamma$, respectively). Therefore,

\begin{equation}
\rho'_2(\ell) =
\lim_{dx \rightarrow 0} {1 \over dx} \int_{I_1 \cap [\beta_1,\beta_2]} d\beta'\ {\sin\beta' \sin(\beta'+\gamma) \over \sin\gamma} +
\lim_{dx \rightarrow 0} {1 \over dx} \int_{I_2 \cap [\beta_1,\beta_2]} d\beta'\ {\sin\beta' \sin(\beta'+\gamma) \over \sin\gamma}.
\end{equation}
Since the intervals $I_1$ and $I_2$ are infinitesimal, they can be made so small that they either are fully contained in $[\beta_1,\beta_2]$ or are disjoint with it. Suppose that $I_1 \subset [\beta_1,\beta_2]$. Then the first integral is, to first order in $dx$, the value of its integrand
when $\beta' = \arcsin\left( {x \over \ell} \sin\gamma \right) - \gamma$ times the width of the interval $I_1$. When $\beta' = \arcsin\left( {x \over \ell} \sin\gamma \right) - \gamma$, then $\sin\beta' \sin(\beta'+\gamma) = \sin\left( \arcsin\left( {x \over \ell} \sin\gamma \right) - \gamma \right) {x \over \ell} \sin\gamma$.
To find the width of the interval $I_1$ we use the expansion

\begin{equation}
\arcsin\left( {x + dx \over \ell} \sin\gamma \right) - \gamma = \arcsin\left( {x \over \ell} \sin\gamma \right) - \gamma + {\sin\gamma \over \sqrt{\ell^2 - x^2 \sin^2\gamma}}\ dx + o(dx).
\end{equation}
Therefore, if $I_1 \subset [\beta_1,\beta_2]$, the first term of expression (64) is

\begin{equation}\label{66}
{x \sin\gamma \over \ell \sqrt{\ell^2 - x^2 \sin^2\gamma}} \sin\left( \arcsin\left( {x \over \ell} \sin\gamma \right) - \gamma \right).
\end{equation}
Likewise, if $I_2 \subset [\beta_1,\beta_2]$, the second term of expression (64) is

\begin{equation}\label{67}
{x \sin\gamma \over \ell \sqrt{\ell^2 - x^2 \sin^2\gamma}} \sin\left( \arcsin\left( {x \over \ell} \sin\gamma \right) + \gamma \right).
\end{equation}
The sine theorem yields $y = x {\sin\beta' \over \sin(\gamma + \beta')}$. The $y$ coordinates which correspond to $\beta' = \arcsin\left( {x \over \ell} \sin\gamma \right) - \gamma$ and to $\beta' = \pi - \arcsin\left( {x \over \ell} \sin\gamma \right) - \gamma$ are

\begin{equation}\label{68}
\begin{array}{cc}
  & y_- = {\sin\left( \arcsin\left( {x \over \ell} \sin\gamma \right) - \gamma \right) \over \sin\gamma}\ \ell \\
  {\rm{and}}&  \\
  & y_+ = {\sin\left( \arcsin\left( {x \over \ell} \sin\gamma \right) + \gamma \right) \over \sin\gamma}\ \ell,
\end{array}
\end{equation}
respectively (note that the functions $y_\pm$ take a complex value when $x \sin\gamma > \ell$, i. e., when $\ell$ is too short to reach the $y$ side). Therefore,

$$
\rho'_2(\ell) = 1_{[y_1,y_2]} (y_-) {x \sin\gamma \over \ell \sqrt{\ell^2 - x^2 \sin^2\gamma}} \sin\left( \arcsin\left( {x \over \ell} \sin\gamma \right) - \gamma \right) +
$$
\begin{equation}\label{69}
1_{[y_1,y_2]} (y_+) {x \sin\gamma \over \ell \sqrt{\ell^2 - x^2 \sin^2\gamma}} \sin\left( \arcsin\left( {x \over \ell} \sin\gamma \right) + \gamma \right),
\end{equation}
where $1_{[y_1,y_2]}$ is the indicator function of the set $[y_1,y_2]$, that is, $1_{[y_1,y_2]} (y) = 1$ if $y \in [y_1,y_2]$ and 0 otherwise.

\begin{figure}[ht]
  \begin{center}
  \includegraphics[width=0.35 \textwidth]{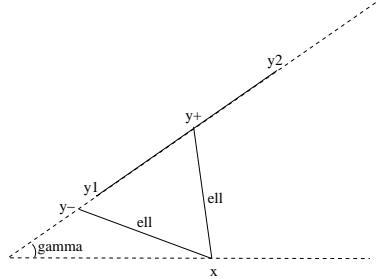}\caption{\label{Conc4.1}In this case $1_{[y_1,y_2]} (y_-) = 0$ and $1_{[y_1,y_2]} (y_+) = 1$}
  \end{center}
\end{figure}

When $\gamma < \pi/2$ there are two segments of length $\ell$ that go from the point of coordinate $x$ to the straight line supporting the coordinate $y$. When $\gamma > \pi/2$ there is only one such segment. Consequently, when $\gamma < \pi/2$, each of the prefactors may be 0 or 1, while if $\gamma > \pi/2$, the first prefactor is always 0.

When $\gamma = 0$ and the distance between the parallel lines is $d$, the remarks previous to formula (\ref{56}) lead to

\begin{equation}
\begin{array}{cc}
  & y_- = x - \sqrt{\ell^2 - d^2}, \\
  & \\
  & y_+ = x + \sqrt{\ell^2 - d^2},
\end{array}
\end{equation}
and

\begin{equation}
\rho'_2(\ell) =  \Big( 1_{[y_1,y_2]} (y_-) + 1_{[y_1,y_2]} (y_+) \Big) {d^2 \over \ell^2 \sqrt{\ell^2 - d^2}}.
\end{equation}

To use the preceding four formulae to compute the density of measure as prescribed in formula (\ref{45}), we note that the angle $\gamma$ is the angle $AOD$ of Figure \ref{PrePajarita} and we reverse the change in notation (\ref{58}):

\begin{equation}
y_1 \rightarrow |OC|,\ \ \ y_2 \rightarrow |O op\left( A + a \vec u \right)|,\ \ \ x \rightarrow |OA - a \vec u|.
\end{equation}

We have proven

{\bf{Proposition 9.2.}}{\it{ The integrand $\rho'_2$ of Proposition 9.1 is given by formulae (\ref{68}) and (\ref{69}). The angle $\gamma$ is the angle $AOD$ of Figure \ref{PrePajarita} and the change of notation ({\ref{58}) is reversed:

\begin{equation}
y_1 \rightarrow |OC|,\ \ \ y_2 \rightarrow |O op\left( A + a \vec u \right)|,\ \ \ x \rightarrow |OA - a \vec u|.
\end{equation}
When $\gamma = 0$ and the distance between the parallel lines is $d$, then the integrand $\rho'_2$ of Proposition 9.1 is given by formulae (\ref{66}) and (\ref{67}). }}

\section{The general case}

The quadrilateral in section 7 had four sides to start with. But in order to compute its cld, its sides $OA$ and $OD$ were conveniently divided into two segments. Together with the sides $VA$ and $VD$ this gave a total of six segments, which gave 15 pairs. Five of the sets of chords joining the segments of each pair were empty, nine were of the kind studied in sections 2-6, and the last set of chords was of a new kind whose cld was derived in sections 7-9. Thus the non empty sets of chords between two of these segments are either the chords between two opposing sides of a convex quadrilateral or the chords between two particular sides of a diagonal pentagon. In this section we prove that such a procedure can be applied to any polygon, that is, the concavity can always be packed in diagonal pentagons, so to speak.

{\bf{Definitions.}}
Let $AB$ and $CD$ be two segments such that their convex cover is a quadrilateral. Then by $Qu(AB, CD)$ we denote the set of chords joining the sides $AB$ and $CD$ of the said convex quadrilateral. Let $AB$ and $CD$ be two segments and $V$ be a point such that $BAVDC$ is a diagonal pentagon whose concave vertex is $V$. Then by $DPent(BA,AVD,DC)$ we denote the set of chords joining the sides $AB$ and $CD$ of the diagonal pentagon $BAVDC$. When $AB$ and $CD$ are sides of a polygon, the segments in $Qu(AB, CD)$ may or may not be fully contained in the polygon. The latter segments constitute the set of chords $Ch(AB, CD)$. By definition, $Ch(AB, CD) \subseteqq Qu(AB, CD)$.

When the polygon is concave there is at least one pair of sides $AB$ and $CD$ for which $\emptyset \subset Ch(AB, CD) \subset Qu(AB, CD)$ holds, where the inclusions are not trivial. Then there is at least one segment $s \in Qu(AB, CD)$ which is a chord. $s$ divides the convex quadrilateral $ABCD$ into two parts (we suppose here that $ABCD$ is the order in which the vertices yield a convex quadrilateral). At least one of the parts must contain a point $P$ which is not in the polygon. Any segment joining the chord $s$ and the point $P$ must cross the perimeter of the polygon at least once. Therefore, in the said half of the quadrilateral lays a portion of the perimeter of the polygon. This portion cannot be just a part of a side, for in that case it would cross at least one of the sides $AB$ and $CD$, which would not be sides themselves. Thus it must contain at least one vertex.

{\bf{Definitions.}} Consider the set which is the union of the points $\{A, D\}$ and the portion of the perimeter of the polygon which is contained in the quadrilateral $ABCD$ and crosses the segment $AD$. We denote its convex cover by $Per(A, D)$. $Per(B, C)$ is defined likewise. By $AB_v$ we denote the intersection of the side $AB$ and the chords in $Ch(AB, CD)$. Informally speaking, it is the part of $AB$ from which some point of $CD$ is visible. $CD_v$ is defined likewise.

{\bf{Proposition 10.1}} {\it{For any two sides $AB$ and $CD$ of a polygon, the side $AB$ can be covered by a finite collection of non overlapping segments with the following property. Let $P_A P_B$ be any of the said segments. Then the set of chords that join $P_A P_B$ and $CD$ is either empty or is the set union or difference of at most three sets of chords of the form $Qu$ or $DPent$ defined above.}}

{\bf{Proof.}} Take any point $P$ in the interior of $AB_v$. Consider the chords $s_1(P)$ and $s_2(P)$ which meet the point $P$ and are tangent to $Per(A, D)$ and $Per(B, C)$, respectively. Almost surely $s_1(P)$ and $s_2(P)$ touch only one vertex $V_1$ of $Per(A, D)$ and one vertex $V_2$ of $Per(B, C)$, respectively. These are the vertices which limit the vision of the side $CD$ from the point $P$. As the point $P$ moves to either side it reaches points $P_A$ and $P_B$ such that the set of vertices touched by $s_1(P)$ and $s_2(P)$ increases by one vertex, because $s_1(P)$ or $s_2(P)$ include some side of $Per(A, D)$ or $Per(B, C)$. Outside the segment $P_{A} P_{B}$ the vision of the side $CD$ is limited by a different pair of vertices. One can proceed in this manner and cover $AB_v$ by a set of disjoint segments such that in each of them the vision of the side $CD$ is limited by a different pair of vertices.

For each such segment $P_A P_B$, there are eight possible orderings of the points $Q_A$, $Q'_A$, $Q_B$, $Q'_B$, and we shall show that for each of them Proposition 10.1 is true.

Case 1. In the case shown in Figure \ref{Conc8}, $Ch(P_A P_B, Q_A Q'_B) = $

\noindent
$Ch(P_A P_B, Q_A Q_B) + Ch(P_A P_B, Q_B Q'_A) + Ch(P_A P_B, Q'_A Q'_B) = $

\noindent
$DPent(P_A P_B, P_B V_2 Q_A, Q_A Q_B) + Qu(P_A P_B, Q_B Q'_A) + DPent(P_B P_A, P_A V_1 Q'_B , Q'_B Q'_A)$.

When $Q_A = Q_B$ and/or $Q'_A = Q'_B$ two/one of the three previous terms disappear and the Proposition 10.1 holds.

\begin{figure}[ht]
  \begin{center}
  \includegraphics[width=0.7 \textwidth]{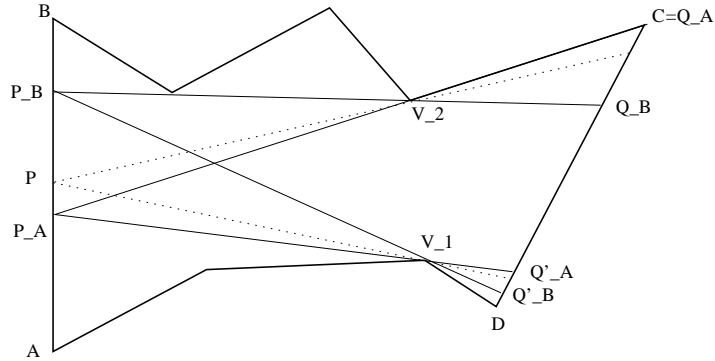}\caption{\label{Conc8}Case 1 of Proposition 10.1}
  \end{center}
\end{figure}

\begin{figure}[ht]
  \begin{center}
  \includegraphics[width=0.7 \textwidth]{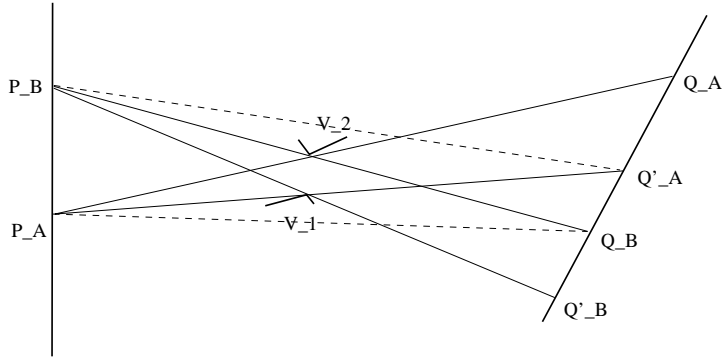}\caption{\label{Conc9}Case 2 of Proposition 10.1}
  \end{center}
\end{figure}

Case 2. In the case shown in Figure \ref{Conc9}, $Ch(P_A P_B, Q_A Q'_B) = $

\noindent
$Ch(P_A P_B, Q_A Q_B) + Ch(P_A P_B, Q'_A Q'_B) - Ch(P_A P_B, Q'_A Q_B) = $

\noindent
$(DPent(P_A P_B, P_B V_2 Q_A, Q_A Q_B) - $ segments that join $P_A P_B$ and $Q_A Q_B$ and go below the vertex $V_1) + (DPent(P_B P_A, P_A V_1 Q'_B, Q'_B Q'_A) - $segments that join $P_A P_B$ and $Q'_A Q'_B$ and go above the vertex $V_2) - (Qu(P_A P_B, Q'_A Q_B) - $segments that join $P_A P_B$ and $Q'_A Q_B$ and go below the vertex $V_1 - $segments that join $P_A P_B$ and $Q'_A Q_B$ and go above the vertex $V_2)$.

But the segments that join $P_A P_B$ and $Q_A Q_B$ and go below the vertex $V_1$ are the segments that join $P_A P_B$ and $Q'_A Q_B$ and go below the vertex $V_1$, and the segments that join $P_A P_B$ and $Q'_A Q'_B$ and go above the vertex $V_2$ are the segments that join $P_A P_B$ and $Q'_A Q_B$ and go above the vertex $V_2$. Therefore, in the preceding equality between sets, the sets that have been named using sentences in English cancel and we obtain

$Ch(P_A P_B, Q_A Q'_B) = $

$Ch(P_A P_B, Q_A Q_B) + Ch(P_A P_B, Q'_A Q'_B) - Ch(P_A P_B, Q'_A Q_B) = $

$DPent(P_A P_B, P_B V_2 Q_A, Q_A Q_B) + DPent(P_B P_A, P_A V_1 Q'_B, Q'_B Q'_A) - Qu(P_A P_B, Q'_A Q_B)$.

Case 3. $Q'_A = Q_B$. This is a limiting case of both case 1 and case 2. Either way we obtain:

$Ch(P_A P_B, Q_A Q'_B) = DPent(P_A P_B, P_B V_2 Q_A, Q_A Q_B) + DPent(P_B P_A, P_A V_1 Q'_B, Q'_B Q'_A)$.

Case 4. $Q'_A = Q_A$. In this case $P_A, V_1, V_2$ and $Q_A = Q'_A$ are aligned. $Ch(P_A P_B, Q_A Q'_B)$ decomposes as in Case 2.

Case 5. $Q'_B = Q_B$. In this case $P_B, V_2, V_1$ and $Q_B = Q'_B$ are aligned. $Ch(P_A P_B, Q_A Q'_B)$ decomposes as in Case 2.

Case 6. $Q_A = Q_B$. In this case $V_2$ is on the side $CD$. $Ch(P_A P_B, Q_A Q'_B)$ decomposes as in Case 1 but has only two terms.

Case 7. $Q'_A = Q'_B$. In this case $V_1$ is on the side $CD$. $Ch(P_A P_B, Q_A Q'_B)$ decomposes as in Case 1 but has only two terms.

Case 8. $Q_A = Q_B$ and $Q'_A = Q'_B$. In this case $V_1$ and $V_2$ are on the side $CD$. $Ch(P_A P_B, Q_A Q'_B)$ decomposes as in Case 1 but has only the $Qu(P_A P_B, Q_B Q'_A)$ term.

When the polygon is not simply connected, there might be holes in the interior of $ABCD$. Each of this holes has a convex cover which we denote by $Per(Hole)$. Then the vision of the side $CD$ from a point $P \in AB$ may consist of more than one segment (up to $h+1$ segments if $h$ holes are caught in the quadrilateral $ABCD$) and instead of two chords $s_1(P)$ and $s_2(P)$ there might now be up to $2 (h+1)$ chords which are tangent to two vertices in $Per(A, D)$ and $Per(B, C)$ and to 2 $h$ vertices in $Per(Hole)$. But each of the segments $Q_A Q'_B$ corresponding to each pair of vertices can still only fall in any of the eight cases considered above. Thus the proposition also holds for not simply connected polygons.

Since the number of vertices in $Per(A, D)$ or $Per(B, C)$ or $Per(Hole)$ is finite, the number of segments $P_A P_B$ is finite. $\ \ \Box$

\section{Concave examples}

\subsection{Concave pentagons}

\begin{figure}[ht]
  \begin{center}
  \includegraphics[width=0.5 \textwidth]{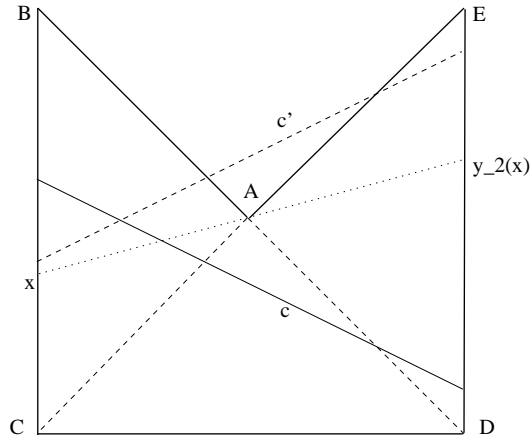}\caption{\label{Sobre}A concave symmetric pentagon}
  \end{center}
\end{figure}

The measure distribution of the shown symmetric diagonal pentagon (Figure \ref{Sobre}) is

$$
2 \rho_{AB, BC} + 2 \rho_{AB, CD} + 2 \rho_{BC, CD} + \rho_{conc\ BC, ED} =
$$

\begin{equation}
2 \rho_{convex\ AB, BC} + 2 \rho_{convex\ AB, CD} + 2 \rho_{convex\ BC, CD} + {1 \over 2} \rho_{convex\ BC, ED},
\end{equation}
because for each chord $c$ joining the sides $BC$ and $ED$ contained in the pentagon there is a mirror image of it, $c'$, joining the sides $BC$ and $ED$ but not contained in the pentagon.

To compute $\rho_{conc\ BC, ED}$ we may also use (70) and (71). If we let $B, C, D$ and $E$ be the vertices of a square of side length 1, then $d=1,\ y_1=0$ and $y_2(x)=1-x$. Then (70) and (71) become

\begin{equation}
\begin{array}{cc}
  & y_- = x - \sqrt{\ell^2 - 1}, \\
  & \\
  & y_+ = x + \sqrt{\ell^2 - 1},
\end{array}
\end{equation}
and

\begin{equation}
\rho'_2(\ell) =  \Big( 1_{[0,1-x]} (y_-) + 1_{[0,1-x]} (y_+) \Big) {1 \over \ell^2 \sqrt{\ell^2 - 1}}.
\end{equation}
Thus

$$
\rho_{conc\ BC, ED}(\ell) = \int_0^1 dx\ \Big( 1_{[0,1-x]}(x - \sqrt{\ell^2 - 1}) + 1_{[0,1-x]}(x + \sqrt{\ell^2 - 1}) \Big) {1 \over \ell^2 \sqrt{\ell^2 - 1}} =
$$

\begin{equation}
\int_0^1 dx\ \Big( 1_{[\sqrt{\ell^2 - 1},(1+\sqrt{\ell^2 - 1})/2]}(x) + 1_{[0,(1-\sqrt{\ell^2 - 1})/2]}(x) \Big) {1 \over \ell^2 \sqrt{\ell^2 - 1}} =
{1-\sqrt{\ell^2 - 1} \over \ell^2 \sqrt{\ell^2 - 1}},
\end{equation}
when $\ell \in [1, \sqrt{2}]$ and 0 otherwise.

For every chord $c$ joining $BC$ and $ED$ and passing below the vertex $A$, there is a symmetrical chord $c'$ of the same length passing above the vertex $A$. Therefore using $\rho_{conc\ BC, ED} = {1 \over 2} \rho_{convex\ BC, ED}$ one obtains the same result as above.

\bigskip

\begin{figure}[ht]
  \begin{center}
  \includegraphics[width=0.6 \textwidth]{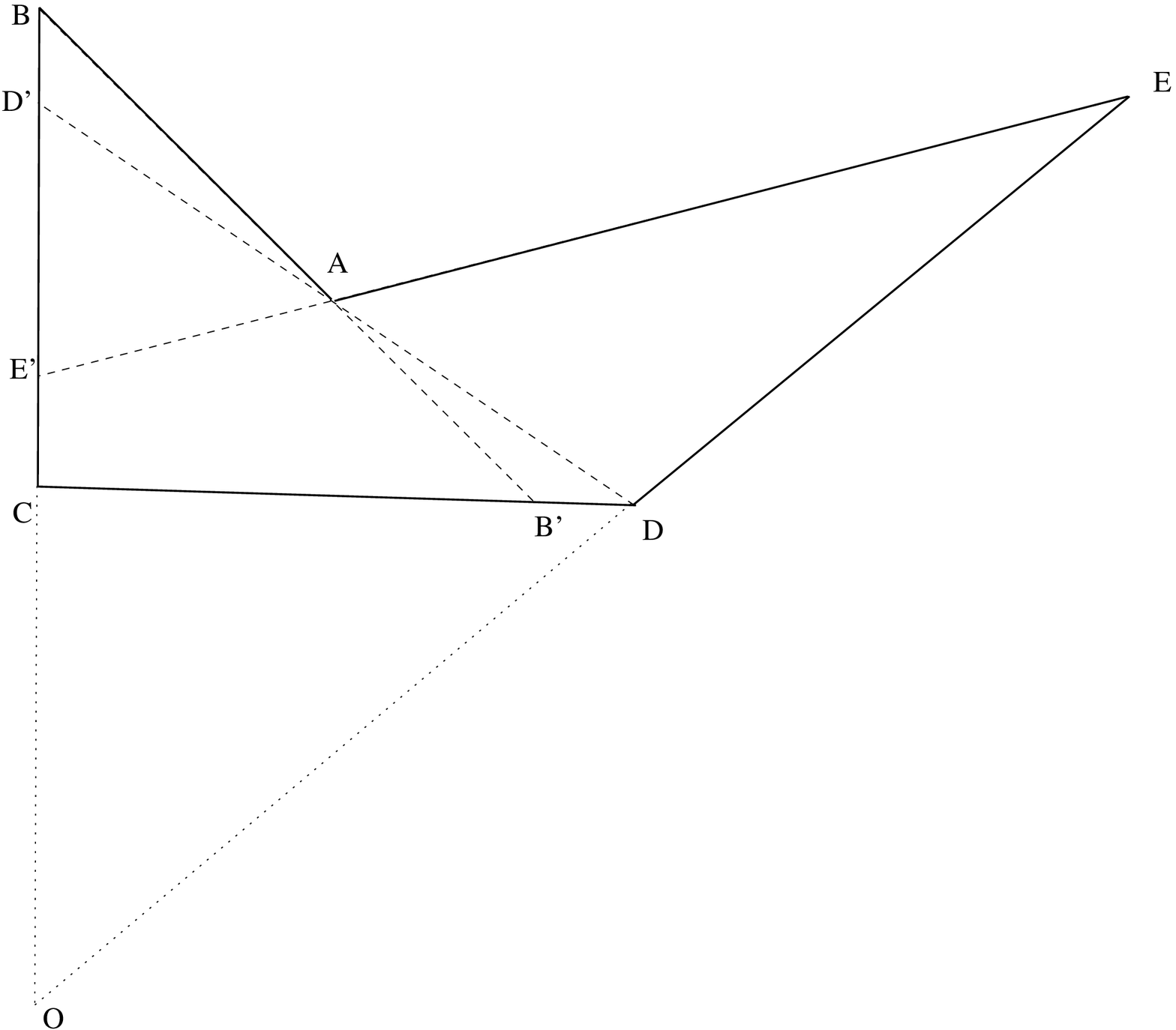}\caption{\label{Sobre2}A concave non symmetric pentagon}
  \end{center}
\end{figure}

Another example is the concave non symmetric, non diagonal pentagon depicted in Figure \ref{Sobre2}. Its measure density has terms which are computed as explained in the part of the article devoted to convex polygons:

$$
\rho_{AB, BC} + \rho_{AB, CB'} + \rho_{BC, CB'} + \rho_{D'C, B'D} + \rho_{E'C, DE} + \rho_{E'C, EA} + \rho_{CD, DE} + \rho_{CD, EA} + \rho_{DE, EA},
$$
and the following two terms, which are computed using (68) and (69):

$$
\rho_{conc\ BD', B'D} + \rho_{conc\ D'E', DE}
$$
When computing $\rho_{conc\ BD', B'D}$, $\gamma = \pi/2\ $. This implies that $y_{\pm} = \pm \sqrt{\ell^2 - x^2}$ and that the first term of formula (69) vanishes (as remarked after formula (69)), leaving

\begin{equation}
\rho'_2(\ell) = 1_{[y_1,y_2]} {x \over \ell^2}.
\end{equation}
If we let the coordinate $x$ be the one along the segment $CD$, with origin in $C$, then

\begin{equation}
\rho_{conc\ BD', B'D} (\ell) = \int_{|CB'|}^{|CD|} dx\ 1_{[CD',Cop(x)]}(\sqrt{\ell^2-x^2}) {x \over \ell^2}.
\end{equation}
In this case the indicator function is 1 when $|CD'| < \sqrt{\ell^2 - x^2} < {y_A x \over x-x_A}$, where $x_A$ and $y_A$ are the Cartesian coordinates of $A$ with respect to a system with origin in $C$ and $x$ and $y$ axes along $CD$ and $CB$, respectively. It is 0 otherwise. The double inequality at the beginning of this paragraph may be multiplied by $x-x_A$ and squared. A double inequality which is cubic in $x$ is then obtained, from which the intervals in which the indicator function is 1 can be solved for and $\rho_{conc\ BD', B'D} (\ell)$ may be computed analytically.

$\rho_{conc\ D'E', DE}$ is a generic case and there is no simplification of functions (67) or (68). The expression for the measure density may be written as follows:

\begin{equation}
\rho_{conc\ D'E', DE} (\ell) = \int_{OE'}^{OD'} dx\ \Big( 1_{[OD,Oop(x)]} \left( y_ - \right) \cdots + 1_{[OD,Oop(x)]} \left( y_+ \right) \cdots \Big),
\end{equation}
where the integrand is a short hand notation for the function (69) showing the extremes of its indicator functions. The angle $\gamma$ is $\widehat{COD}$.

\subsection{The five-pointed star}

\begin{figure}[ht]
  \begin{center}
  \includegraphics[width=0.45 \textwidth]{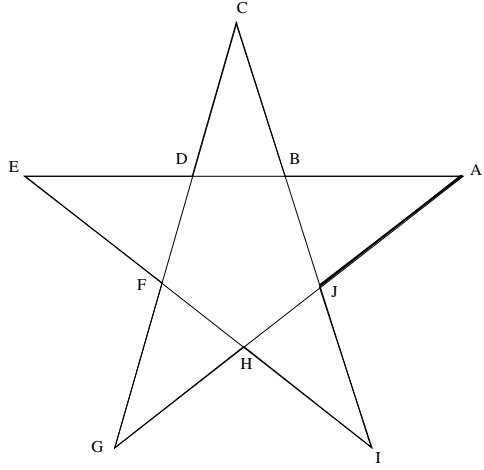}\caption{\label{Pentagrama2}The five-pointed star}
  \end{center}
\end{figure}

In Figure \ref{Pentagrama2}, the polar coordinates of the outer vertices of the five-pointed star are, starting from A, $(1, \pi/10),\ (1, \pi/2),\ (1, 9 \pi/10),\ (1, 13 \pi/10)$ and $(1, 17 \pi/10)$. The polar coordinates of the inner vertices of the five-pointed star are, starting from B, $(r, 3 \pi/10),\ (r, 7 \pi/10),$ $(r, 11 \pi/10),$ $(r, 3 \pi/2)$ and $(r, 19 \pi/10)$, where $r = {\sin (\pi/10) \over \sin (7 \pi/10)}$. Due to the symmetry of the star, the set of chords which touch the side $JA$ is congruent to the set of chords which touch any other side. If we multiplied the said set by 10, we would count each chord twice, because all chords (up to a set of measure 0) would touch two sides. Therefore the density of measure for the star is five times the the density of measure of the chords which touch the side $JA$. In its turn,

\begin{equation}
\rho_{JA} = \rho_{JA,AB} + \rho_{conc\ AJ,CD} + \rho_{JA,DE} + \rho_{JA,EF} + \rho_{JA,FG}.
\end{equation}
A rotation of $3 (2 \pi/5)$ transforms the pair of sides $JA, DE$ into the pair of sides $FG, JA$. Therefore:

\begin{equation}
\rho_{JA} =  \rho_{JA,AB} + \rho_{conc\ AJ,CD} + 2 \rho_{JA,DE} + \rho_{JA,EF},
\end{equation}
where

\begin{equation}
\rho_{conc\ AJ,CD} (\ell) = \int_{|J-G|}^{|A-G|}dx\ \Big( 1_{[GD,Gop(x)]} \left( y_ - \right) \cdots + 1_{[GD,Gop(x)]} \left( y_+ \right) \cdots \Big).
\end{equation}
The functions $y_{\pm}$ are given by formula (68) with $\gamma = \pi/5$ and the dots stand for

\noindent
${x \sin(\pi/5) \over \ell \sqrt{\ell^2 - x^2 \sin^2(\pi/5)}} \sin\left( \arcsin\left( {x \over \ell} \sin{\pi \over 5} \right) \pm {\pi \over 5} \right)$, where the sign of the last $\pi \over 5$ matches the subindex of the argument of the corresponding indicator function.

The above function is expensive to compute. The author has approximated the Riemann integral above by a Riemann sum with an interval of 0.00001. The computation of the graph that follows took 1,100 seconds in his personal computer, using the software Mathematica. Even after 1,100 one can still see some spikes due to numerical (not real) oscillations. The numerical nature of the spikes is revealed by the fact that they spread outside the bound $[app_{-}(n)(\ell), app_{+}(n)(\ell)]$ (see section \ref{Converging bounds}). These spikes arise because for a given width of the Riemann sum the error is a rapidly varying function of $\ell$.

It is numerically better to use the method given in section \ref{Converging bounds}. When $n = 1000$, the graphs of $app_{-}(n)(\ell)$ and $app_{+}(n)(\ell)$ for $\rho_{conc\ AJ,CD}$ are indistinguishable from each other. It takes 180 seconds to compute each of them, on the same computer as before. These bounds are exact in the sense that they are computed in terms of measure densities for convex polygons, which in turn we had computed in terms of elementary functions.

We plot in Fig. \ref{Superposition} their average and the superposition with the graph of $\rho_{conc\ AJ,CD}(\ell)$ computed using the Riemann sum. Both ways of computing $\rho_{conc\ AJ,CD}$ are worked out in the file Concave.nb in \cite{RGPN}.

\begin{figure}[ht]
  \begin{center}
  \includegraphics[width=0.4 \textwidth]{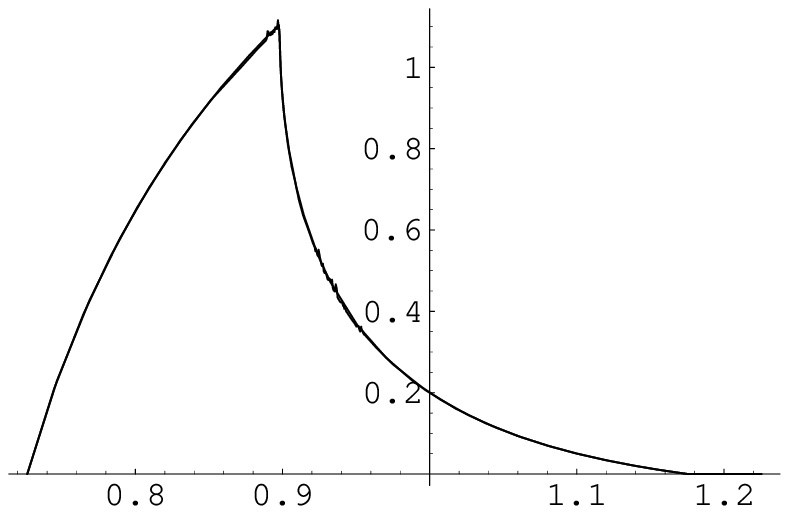}\caption{\label{Superposition}Superposition of the graphs of $\rho_{conc\ AJ,CD}(\ell)$ computed using (67) and (68) (which yields some small spikes) and using (39)}
  \includegraphics[width=0.3 \textwidth]{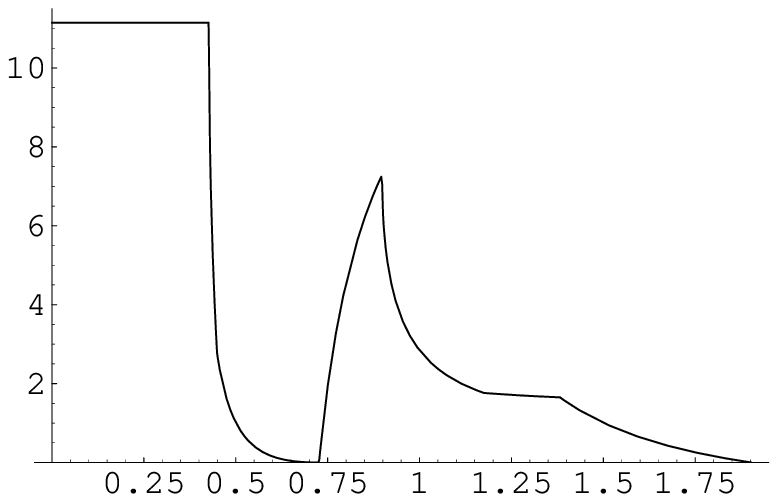}\caption{\label{PentagramaGraf}Measure density of the length of chords for the five-pointed star}
  \end{center}
\end{figure}

In Figure \ref{PentagramaGraf} the graph of the measure density (i. e., $5 \rho_{JA}(\ell)$) is plotted.


\section{Remark}
We notice that

$$
\int_{x_1}^{x_2} dx\ 1_{[y_1, y_2]}(y_{\pm}(x, \ell)) {x \sin\gamma \over \ell \sqrt{\ell^2 - x^2 \sin^2\gamma}} \sin\left( \arcsin\left( {x \over \ell} \sin\gamma \right) \pm \gamma \right) =
$$

\begin{equation}
\int_{[x_1, x_2] \cap \{x| y_\pm(x, \ell) \in [y_1, y_2]\}} dx\ {x \sin\gamma \over \ell \sqrt{\ell^2 - x^2 \sin^2\gamma}} \sin\left( \arcsin\left( {x \over \ell} \sin\gamma \right) \pm \gamma \right)
\end{equation}
and that the primitives of ${x \sin\gamma \over \ell \sqrt{\ell^2 - x^2 \sin^2\gamma}} \sin\left( \arcsin\left( {x \over \ell} \sin\gamma \right) \pm \gamma \right)$ are

\begin{equation}
-{\sin\gamma \over \ell} \left( {x \over 2 \ell} \left( x \sin\gamma \mp \cot\gamma \sqrt{\ell^2 - x^2 \sin^2\gamma} \right) + {\ell \cos\gamma \over 2 \sin^2\gamma} \arctan {\sqrt{\ell^2 - x^2 \sin^2\gamma} \over x \sin\gamma} \right).
\end{equation}

The difficulty in computing the above integral then, reduces to the computation of the integration limits. These limits are quite laborious to find, and have different analytical forms which depend on $\ell$ and on the geometry of the pair of segments. Still, it is possible in principle to obtain results in the manner of H. S. Harutyunyan, V. K. Ohanyan and A. Gasparyan \cite{Harutyunyan2011,Gasparyan2013}. While this is less practical as a general purpose procedure to find distributions of chords, it may be worthwhile doing if one wishes to compute distributions of chords for a certain class of polygons repeatedly.

\bigskip
{\large\bf{Appendix: Opposite function}}

\begin{figure}[h]
  \begin{center}
  \includegraphics[width=0.5 \textwidth]{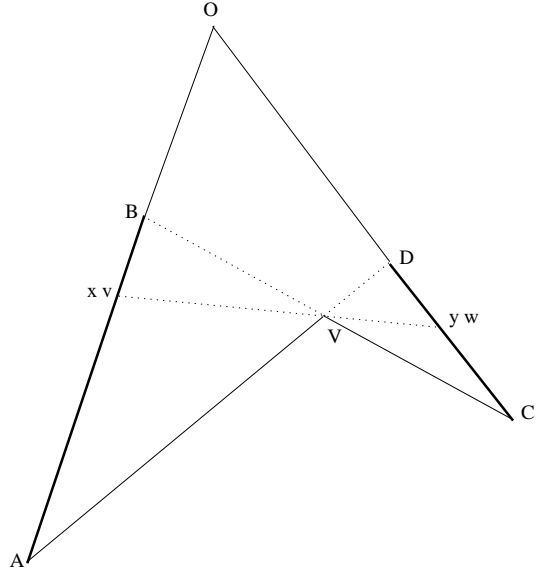}\caption{\label{Op}The opposite function}
  \end{center}
\end{figure}

We define the unit vectors

\begin{equation}
v \equiv {A - B \over |A - B|}\ \ \ \ {\rm{and}}\ \ \ \ w \equiv {C - D \over |C - D|}
\end{equation}
and let $x$ and $y$ be the distances to $O$ along the straight lines $OA$ and $OC$, respectively.

We want to obtain the coordinates of $O + y w$ when the coordinates of $A, B, C, D$ and $O + x\ v$ are known. Two ways to obtain them are the following: Solve $y$ from the equation
\begin{equation}
(y\ w + O) - (x\ v + O) \propto V - (x\ v + O)
\end{equation}
or notice that the parallelograms which the vector $O + x\ v - V$ spans with the vectors $D - V$ and $D - (O + y\ w)$ have equal areas, because their difference is parallel to $O + x v - V$. By either way one arrives at

\begin{equation}
(O + x\ v - V) \times (D - V) = (O + x\ v - V) \times (D - O - y\ w),
\end{equation}
from which $y$ can be solved for:

\begin{equation}
op(x) \equiv y = |OD| - {(O + x\ v - V) \times (D - V) \over (O + x\ v - V) \times w}.
\end{equation}

When the segments $AB$ and $CD$ are parallel, $\gamma \equiv \widehat{BOD} = 0$, $w = v$ and $O$ lies at infinity. According to the previous scheme for $\gamma \neq 0$, the $x$ and $y$ coordinates would take infinite values at the segments. By convention, in this case, we take the origin of the $x$ coordinates at $B$ and the origin of the $y$ coordinates at the closest point to $B$ which lies on the straight line which contains the segment $CD$. As in the $\gamma \neq 0$ case the chords which lie in the polygon under consideration are the ones which cross the triangle $BVD$. By the same argument as in the $\gamma \neq 0$ case and with the said convention about the origin of coordinates $x$ and $y$,

\begin{equation}
op(x) \equiv y = (D-B).v - {(x\ v - (V - B)) \times (D - V) \over (x\ v - (V - B)) \times v}.
\end{equation}

\vspace{2cm}
{\large\bf{Acknowledgments}}

This work is part of the research project entitled ``Dynamical Analysis, Advanced Orbit Propagation and Simulation of Complex Space Systems" (ESP2013-41634-P) supported by the Spanish Ministry of Economy and Competitiveness. The author thanks the Spanish Government for its financial support.


\begin{thebibliography}{99}

\bibitem{Mazzolo2003}
{A. Mazzolo, B. Rooesslinger} and C.~Diop.
\newblock {On the properties of the chord length distribution, from integral
  geometry to reactor physics}.
\newblock {\em Annals of Nuclear Energy}, 30:1391--1400, 2003.

\bibitem{Aharonyan2005}
N.~G. Aharonyan and V.~K. Ohanyan.
\newblock {Chord length distribution functions for polygons}.
\newblock {\em J. Contemp. Math. Anal.}, 40.

\bibitem{Ambartzumian1990}
R.~Ambartzumian.
\newblock {\em {Factorization Calculus and Geometric Probability}}.
\newblock Cambridge University Press, 1990.

\bibitem{Baesel2008}
U.~Baesel.
\newblock {\em {Geometrische Wahrscheinlichkeiten für nichtkonvexe
  Testelemente}}.
\newblock 2008.

\bibitem{Baesel2014}
U.~Baesel.
\newblock {Random chords and point distances in regular polygons}.
\newblock {\em Acta Math. Univ. Comenianae}, pages 1--18, 2014.

\bibitem{Boettcher2012}
R.~Böttcher.
\newblock {Measure of long segments intersecting both sides of a kite as a
  basis for arbitrary pairs of segments}.
\newblock {\em Seminarberichte aus der Fakultät für Mathematik und Informatik
  der FernUniversität in Hagen}, 85:127--150, 2012.

\bibitem{Lang2001}
{C. Lang, J. Ohser} and R.~Hilfer.
\newblock {On the analysis of spatial binary images}.
\newblock {\em Journal of Microscopy}, 203:303--313, 2001.

\bibitem{Cauchy1850}
A.~Cauchy.
\newblock {Mémoire sur la rectification des courbes et la quadrature des
  surfaces courbes (1850)}.
\newblock {\em Oeuvres complètes, vol. 2}.

\bibitem{PrivateCommunication}
S.~Ciccariello.
\newblock {Private communication}.
\newblock 2015.

\bibitem{Ciccariello2009}
S.~Ciccariello.
\newblock {The Correlation Function of Plane Polygons}.
\newblock {\em J. Math. Phys.}, 50:103527, 2009.

\bibitem{Ciccariello2010}
S.~Ciccariello.
\newblock {The isotropic correlation function of plane figures: the triangle
  case}.
\newblock {\em Journal of Physics: Conference Series}, 247:012014, 2010.

\bibitem{Clark2002}
M.~Clark.
\newblock {\em {Paradoxes from A to Z}}.
\newblock Routledge, London, 2002.

\bibitem{Coleman1969}
R.~Coleman.
\newblock Random paths through convex bodies.
\newblock {\em Journal of Applied Probability}, 6(2):430--441, 1969.

\bibitem{Czuber1884}
E.~Czuber.
\newblock {Zur Theorie der geometrischen Wahrscheinlichkeiten}.
\newblock {\em Sitzungerber. Akad. Wiss. Abl. 2}, 90:719--742, 1884.

\bibitem{Deltheil1919}
R.~Deltheil.
\newblock {Sur la th\'eorie des probabilit\'es g\'eom\'etriques}.
\newblock {\em Ann. Fac. Sci. Univ. Toulouse (3)}, 11:1--65, 1919.

\bibitem{Duma2009}
A.~Duma and S.~Rizzo.
\newblock {Chord length distribution functions for an arbitrary triangle}.
\newblock {\em Supplemento ai rendiconti del Circolo matematico di Palermo},
  81:141--157, 2009.

\bibitem{Duma2011}
A.~Duma and S.~Rizzo.
\newblock {La funzione di distribuzione di una corda in un trapezio
  rettangolo}.
\newblock {\em Supplemento ai rendiconti del Circolo matematico di Palermo},
  83:147--160, 2011.

\bibitem{Ciccariello1981}
S.~Ciccariello et~al.
\newblock {Correlation functions of amorphous multiphase systems}.
\newblock {\em Phys. Rev. B}, 23(12):6474--6485, 1981.

\bibitem{RGPN}
R.~Garc\'ia-Pelayo.
\newblock {https://sites.google.com/site/ricardogarciapelayo/}.

\bibitem{Gasparyan2013}
A.~Gasparyan and V.~K. Ohanyan.
\newblock {Recognition of triangles by covariogram}.
\newblock {\em Stochastic and Integral Geometry}, 48:110--122, 2013.

\bibitem{Gille2000}
W.~Gille.
\newblock {Chord length distributions and small-angle scattering}.
\newblock {\em Eur. Phys. J. B}, 17:371--383, 2000.

\bibitem{Gille2001a}
W.~Gille.
\newblock {Chord length distribution density of an infinitely long circular
  hollow cylinder}.
\newblock {\em Mathematical and computer modelling}, 34:423--431, 2001.

\bibitem{Gille2001b}
W.~Gille.
\newblock {The small angle scattering correlation function of two infinitely
  long parallel circular cylinders}.
\newblock {\em Computational Materials Science}, 20:181--195, 2001.

\bibitem{Gille2009}
W.~Gille.
\newblock {Geometric figures with the same chord length probability density
  function: Six examples}.
\newblock {\em Powder Technology}, 192:85--91, 2009.

\bibitem{Gruy2008}
F.~Gruy and S.~Jacquier.
\newblock {The chord length distribution of a two-sphere aggregate}.
\newblock {\em Computational Materials Science}, 44:218--223, 2008.

\bibitem{Gruy2014}
F.~Gruy and S.-H. Suh.
\newblock {The chord length distribution of a dumbbell shaped aggregate:
  Analytical expression}.
\newblock {\em Powder Technology}, 253:207--215, 2014.

\bibitem{Harutyunyan2007}
H.~S. Harutyunyan.
\newblock {\em Uchenye Zapiski Yerevan State Univ.}, 1:17--24, 2007.

\bibitem{Harutyunyan2009}
H.~S. Harutyunyan and V.~K. Ohanyan.
\newblock {The chord length distribution function for regular polygons}.
\newblock {\em Adv. Appl. Prob.}, 41:358--366, 2009.

\bibitem{Harutyunyan2011}
H.~S. Harutyunyan and V.~K. Ohanyan.
\newblock {The chord length distribution function for convex polygons}.
\newblock {\em Sutra: International Journal of Mathematical Science Education},
  4:1--15, 2011.

\bibitem{Kendall1963}
M.~G. Kendall and P.~A.~P. Moran.
\newblock {\em {Geometrical Probability}}.
\newblock Charles Griffin \& Company Limited, London, 1963.

\bibitem{Kingman1965}
J.~F.~C. Kingman.
\newblock {Mean free paths in a convex reflecting region}.
\newblock {\em Journal of Applied Probability}, 2:162--168, 1965.

\bibitem{Kingman1969}
J.~F.~C. Kingman.
\newblock {Random Secants of a Convex Body}.
\newblock {\em Journal of Applied Probability}, 6:660--672, 1969.

\bibitem{Mallows1970}
C.~L. Mallows and J.~M.~C. Clark.
\newblock {Linear-Intercept Distributions Do Not Characterize Plane Sets}.
\newblock {\em Journal of Applied Probability}, 7(1):240--244, 1970.

\bibitem{Piefke1978}
F.~Piefke.
\newblock {Beziehungen zwischen der Sehnenl$\rm{\ddot{a}}$ngenverteilung und
  der Verteilung des Abstandes zweier zuf$\rm{\ddot{a}}$lliger Punkte im
  Eik$\rm{\ddot{o}}$rper}.
\newblock {\em Z. Wahrscheinlichkeitstheorie verw. Gebiete}, 43:129--134, 1978.

\bibitem{Piefke1979}
F.~Piefke.
\newblock {Chord length distribution of the ellipse}.
\newblock {\em Lithuanian Mathematical Journal}, 19(3):325--333, 1979.

\bibitem{Ren2014}
D.~Ren.
\newblock {Random chord distributions and containment functions}.
\newblock {\em Adv. Applied Math.}, 58:1--20, 2014.

\bibitem{Barrilla2011}
\rm{D. Barrilla, A. Duma} and A.~Puglisi.
\newblock {The distribution function of a chord in a non-convex poligon}.
\newblock {\em Supplemento ai rendiconti del Circolo matematico di Palermo},
  83:23--40, 2011.

\bibitem{Santalo1976}
L.~A. Santal\'o.
\newblock {\em {Integral Geometry and Geometric Probability}}.
\newblock Addison-Wesley, Reading, MA, 1976.

\bibitem{Sorrenti2012b}
L.~Sorrenti.
\newblock {Chord length distribution functions for an isosceles trapezium}.
\newblock {\em General Mathematics}, 20:9--24, 2012.

\bibitem{Sorrenti2012}
L.~Sorrenti.
\newblock {The chord length distribution function for a non convex polygon}.
\newblock {\em Supplemento ai rendiconti del Circolo matematico di Palermo},
  84:241--250, 2012.

\bibitem{Sulanke1961}
R.~Sulanke.
\newblock {Die Verteilung des Sehnenl$\rm{\ddot{a}}$ngen an ebenen und
  r$\rm{\ddot{a}}$umlichen Figuren}.
\newblock {\em Mathematische Nachrichten}, 23:51--74, 1961.

\bibitem{Vlasov2007}
A.~Y. Vlasov.
\newblock {Signed chord length distribution}.
\newblock {\em e-print arXiv:0711.4734 [math-ph]}.

\bibitem{Vlasov2011}
A.~Y. Vlasov.
\newblock {Extension of Dirac's chord method to the case of a nonconvex set by
  the use of quasi-probability distributions}.
\newblock {\em Journal of Mathematical Physics}, 52:053516, 2011.

\end{thebibliography}

\end{document}